\documentclass[aap,MSNbibl,seceqn,citesort,dvips]{arximspdf}
\usepackage{graphicx}

\doi{10.1214/11-AAP834} 
\volume{22}
\issue{6}
\pubyear{2012}
\firstpage{2357}
\lastpage{2387}

\makeatletter

\newcommand{\cal}{\mathcal}

\newcommand{\calA}{\mathcal{A}}
\newcommand{\calB}{\mathcal{B}}
\newcommand{\calD}{\mathcal{D}}
\newcommand{\calF}{\mathcal{F}}
\newcommand{\calM}{\mathcal{M}}
\newcommand{\calT}{\mathcal{T}}
\newcommand{\calW}{\mathcal{W}}
\newcommand{\calY}{\mathcal{Y}}
\newcommand{\al}{\alpha}
\newcommand{\tha}{\theta}
\newcommand{\ka}{\kappa}
\newcommand{\dUp}{\partial\Upsilon}
\newcommand{\Up}{\Upsilon}
\newcommand{\Om}{\Omega}
\newcommand{\E}{\mathbb{E}}
\newcommand{\N}{\mathbb{N}}
\newcommand{\Pb}{\mathbb{P}}
\newcommand{\R}{\mathbb{R}}
\newcommand{\Z}{\mathbb{Z}}
\newcommand{\bfi}{\mathbf{i}}
\newcommand{\bfj}{\mathbf{j}}
\newcommand{\bfk}{\mathbf{k}}
\newcommand{\bfu}{\mathbf{u}}
\newcommand{\bfv}{\mathbf{v}}
\newcommand{\bfx}{\mathbf{x}}
\newcommand{\Mbar}{\bar{M}}
\newcommand{\half}{{\frac{1}{2}}}

\newcommand{\convergeprob}{\stackrel{\mathbb{P}}\longrightarrow}

\newcommand{\convergeas}{\stackrel{\mathrm{a.s.}}{\longrightarrow}}
\newcommand{\eqdist}{\stackrel{\mathrm{d}}{=}}

\newtheorem{theorem}{Theorem}[section]
\newtheorem{lem}{Lemma}[section]

\newproclaim{assp}{Assumption}[section]
\newproclaim{definition}{Definition}[section]
\newproclaim{remark}{Remark}[section]

\makeatother

\begin{document}
\begin{frontmatter}

\title{A class of multifractal processes constructed using an embedded branching process}
\runtitle{Multifractal process with embedded branching process}

\begin{aug}
\author[A]{\fnms{Geoffrey} \snm{Decrouez}\ead[label=e1]{dgg@unimelb.edu.au}}
\and
\author[A]{\fnms{Owen Dafydd} \snm{Jones}\corref{}\ead[label=e2]{odjones@unimelb.edu.au}}
\runauthor{G. Decrouez and O. D. Jones}
\affiliation{University of Melbourne}
\address[A]{Department of Mathematics\\
\quad and Statistics\\
University of Melbourne\\
Parkville VIC 3010\\
Australia\\
\printead{e1}\\
\phantom{E-mail: }\printead*{e2}} 
\end{aug}

\received{\smonth{10} \syear{2010}}
\revised{\smonth{10} \syear{2011}}

%
\begin{abstract}
We present a new class of multifractal process on $\mathbb{R}$,
constructed using an embedded branching process. The construction makes
use of known results on multitype branching random walks, and along the
way constructs cascade measures on the boundaries of multitype
Galton--Watson trees. Our class of processes includes Brownian motion
subjected to a continuous multifractal time-change.

In addition, if we observe our process at a fixed spatial resolution,
then we can obtain a finite Markov representation of it, which we can
use for on-line simulation. That is, given only the Markov
representation at step $n$, we can generate step $n+1$ in $O(\log n)$
operations. Detailed pseudo-code for this algorithm is provided.
\end{abstract}

%

\end{frontmatter}

\section{Introduction}

Information about the local fluctuations of a process $X$ can be
obtained using the local exponent $h_X(t)$, defined as~\cite{R03}
\[
h_X(t):= \liminf_{\varepsilon\rightarrow0} \frac{1}{\log\varepsilon
}\log
\sup_{|u-t|<\varepsilon} |X(u)-X(t)|.
\]
When $h_X(t)$ is constant all along the sample path with probability
$1$, $X$ is said to be monofractal.
In contrast, we can consider a class of processes whose exponents
behave erratically with time: each interval of positive length exhibits
a range of different exponents.
For such processes, it is, in practice, impossible to estimate $h_X(t)$
for all $t$, due to the finite precision of the data.
Instead, we use the Hausdorff spectrum $D(h)$, a global description of
its local fluctuations.
$D(h)$ is defined as the Hausdorff dimension of the set of points with
a given exponent $h$.
For monofractal processes, $D(h)$ degenerates to a single point at some
$h=H$ [so $D(H)=1$, and the convention is to set $D(h)=-\infty$ for
$h\neq H$].
When the spectrum is nontrivial for a range of values of $h$, the
process is said to be multifractal.

The term multifractal is also well defined for measures.
Let $B(x,r)$ be a ball centered at $x\in\R^n$ with radius $r$.
The local dimension of a finite measure $\mu$ at $x\in\R^n$ is
defined as
\[
\mathrm{dim}_{\mathrm{loc}}\mu(x) = \lim _{r\to0}\frac
{\log\mu
(B(x,r))}{\log r}.
\]
The Hausdorff spectrum $D(\alpha)$ of a measure at scale $\al$ is then
defined as the Hausdorff dimension of the set of points with a given
local dimension $\alpha$.
Measures for which the Hausdorff spectrum does not degenerate to a
point are called multifractal measures.
Constructions of multifractal measures date back to the $m$-ary
cascades of Mandelbrot~\cite{M74}, and the multifractal spectrum of
such measures can be found in, for example,~\cite{R03}.

A positive nondecreasing multifractal process can be obtained by
integrating a multifractal measure.
Other processes with nontrivial multifractal structure can be obtained
by using the integrated measure as a multifractal time change, applied
to monofractal processes such as fractional Brownian motion.
This is the basis of models such as infinitely divisible cascades
\cite{BM02,BM03,CRA05}.

Multifractals have a wide range of applications.
For example, the rich structure of network traffic exhibits
multifractal patterns~\cite{ABFRV02}, as does the stock market
\mbox{\cite{M97,M99}}.
Other applications include turbulence~\cite{SM88}, seismology \cite
{H01,TLM04} and imaging~\cite{RRCB00}, to cite but a few.

On-line simulation of multifractal processes is in general difficult,
because their correlations typically decay slowly, meaning that to
simulate $X(n+1)$ one requires $X(1), \ldots, X(n)$.
This is the same problem faced when simulating fractional Brownian
motion, where to simulate $X(n+1)$ one needs the whole covariance
matrix of $X(1), \ldots, X(n+1)$.
Some simple monofractal processes avoid this problem, for example,
$\alpha$-stable or $M/G/\infty$ processes~\cite{C84}, but it remains a
real problem to find flexible multifractal models that can be quickly simulated.

We propose a new class of multifractal processes, called Multifractal
Embedded Branching Process (MEBP) processes, which can be efficiently
simulated on-line.
MEBP are defined using the crossing tree, an ad-hoc space--time
description of the process, and are such that the spatial component of
their crossing tree is a Galton--Watson branching process.
For any suitable branching process, there is a family of
processes---identical up to a continuous time change---for which the
spatial component of the crossing tree coincides with the branching process.
We identify one of these as the Canonical Embedded Branching Process
(CEBP), and then construct MEBP from it using a multifractal time change.
To allow on-line simulation of the process, the time change is
constructed from a multiplicative cascade on the crossing tree.
The simulation algorithm presents nice features since it only requires
$O(n\log n)$ operations and $O(\log n)$ storage to generate $n$ steps,
and can generate a new step on demand.\vadjust{\goodbreak}

To construct the time change we use here, we start by constructing a
multiplicative cascade on a \textit{multitype} Galton--Watson tree.
The cascade defines a measure on the boundary of the tree, whose
existence follows from known results for multitype branching random walks.
(See, e.g.,~\cite{L00} for the single-type case.)
To map the cascade measure onto $\R_+$, we use the so called
``branching measure'' on the tree, in contrast to the way this is
usually done, using a ``splitting measure.''
See Section~\ref{CEBPMEBP} for details and further background.

The MEBP processes constructed here include a couple of special cases
of interest.
We can represent Brownian motion as a CEBP, thus MEBP processes include
a subclass of multifractal time changed Brownian motions.
Such models are of particular interest in finance~\cite{M97,M99}.
In the special case when the number of subcrossings is constant and
equal to two (for the definition see Section~\ref{secEBP}), the CEBP
degenerates to a straight line, and the time change is just the
well-known binary cascade (see, e.g.,~\cite{B99,KP76,Mo96} and
references therein).

Although we do show that MEBP possess a form of discrete multifractal
scaling [see the discussion following equation (\ref{MFCmulteqn})], the
multifractal nature of MEBP processes is not studied in this paper.
We refer the reader to a coming paper for a full study of the
multifactal spectrum of MEBP~\cite{DHJ}.
In particular, it can be shown that CEBP processes are monofractal, and
that the multifractal formalism holds for MEBP processes, with a
nontrivial spectrum.
The monofractal nature of CEBP processes, together with an upper bound
of the spectrum of MEBP, was derived in the Ph.D. thesis of the first
author~\cite{Dec09}.

The paper is organized as follows.
First we recall the definition of the crossing tree and then construct
the CEBP process.
We then construct MEBP processes and give conditions for continuity.
Finally we provide an efficient on-line algorithm for simulating MEBP processes.
An implementation of the algorithm is available from the second
author's website~\cite{J}.

\section{CEBP and the crossing tree}
\label{secEBP}

Let $X\dvtx\R^+\rightarrow\R$ be a continuous process, with $X(0)=0$.
For $n\in\Z$ we define level $n$ passage times $T_{k}^{n}$ by putting
$T_0^n = 0$ and
\[
T_{k+1}^{n}= \inf\{t>T_{k}^{n}|X(t)\in2^{n}\Z, X(t)\not=
X(T_{k}^{n}) \}.
\]
The $k$th level $n$ (equivalently scale $2^n$) crossing $C_k^n$ is the
sample path from $T_{k-1}^{n}$ to $T_{k}^{n}$.
\[
C_k^n:= \{ (t,X(t))| T_{k-1}^n\leq t < T_k^n\}.
\]

\begin{figure}

\includegraphics{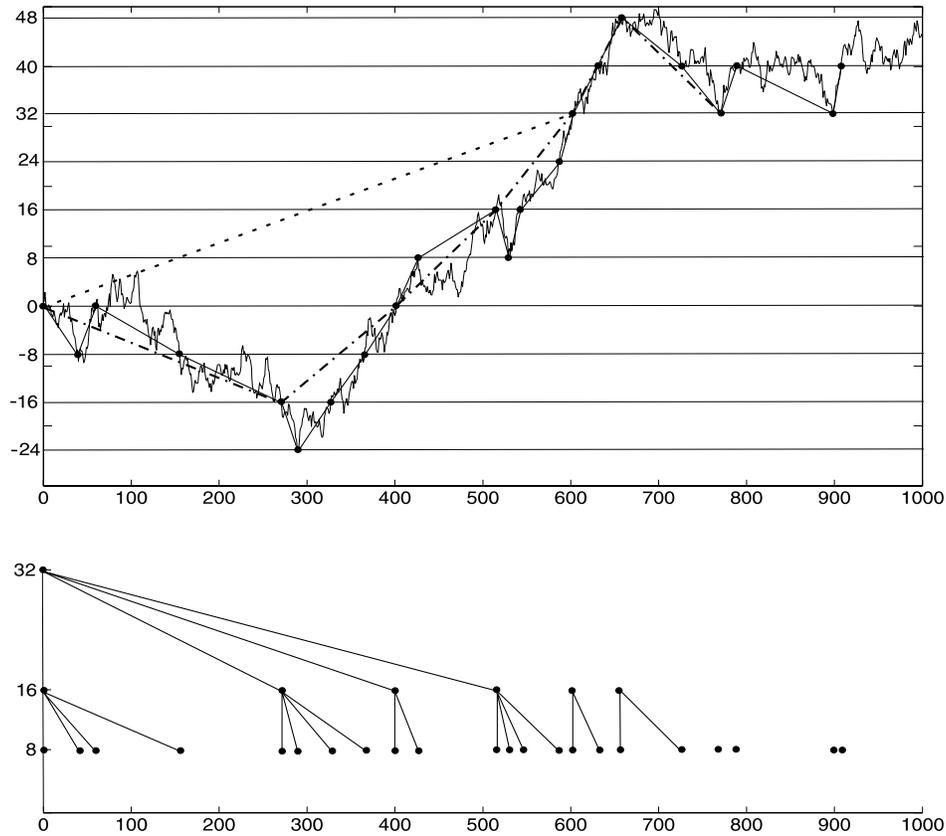}

\caption{A section of sample path and levels 3, 4 and 5 of its crossing
tree. In the top frame we have joined the points $T^n_k$ at each level,
and in the bottom frame we have identified the $k$th level $n$
crossing with the point $(2^n, T^n_{k-1})$ and linked each crossing to
its subcrossings.}%
\label{crossingtree}
\end{figure}

When passing from a coarse scale to a finer one, we decompose each
level $n$ crossing into a sequence of level $n-1$ crossings.
To define the crossing tree, we associate nodes\vadjust{\goodbreak} with crossings, and the
children of a node are its subcrossings.
The crossing tree is illustrated in Figure~\ref{crossingtree}, where
the level 3, 4 and 5 crossings of a given sample path are shown.

The crossing tree is an efficient way of representing a self-similar
signal, and can also be used for inference.
In~\cite{JS04} the crossing tree is used to test for self-similarity
and to obtain an asymptotically consistent estimator of the Hurst index
of a self-similar process with stationary increments, and in \cite
{JS05} it is used to test for stationarity.

In addition to indexing crossings be their level and position within
each level, we will also use a tree indexing scheme.
Let $\varnothing$ be the root of the tree, representing the first level 0
crossing.
The first generation of children (which are level $-1$ crossings, of
size $1/2$) are labeled by $i$, $1\leq i \leq Z_\varnothing$, where
$Z_\varnothing$ is the number of children of~$\varnothing$.
The second generation (which are level $-2$ crossings, of size $1/4$)
are then labeled $ij$, $1\leq j\leq Z_i$, where $Z_i$ is the number of
children of $i$.
More generally, a node is an element of $U = \bigcup_{n\geq0}\N^{n}$ and
a branch is a couple $(\bfu,\bfu j)$ where $\bfu\in U$ and $j\in\N$.
The length of a node $\bfi=i_1, \ldots, i_n$ is $|\bfi|=n$, and the
$k$th element is $\bfi[k] = i_k$.
If $|\bfi|>n$, $\bfi|_n$ is the curtailment of $\bfi$ after $n$ terms.
Conventionally $|\varnothing|=0$, and $\bfi|_0=\varnothing$.
A tree $\Up$ is a set of nodes, that is, a subset of $U$, such that:
\begin{itemize}
\item$\varnothing\in\Up$;
\item if a node $\bfi$ belongs to the tree, then every ancestor node
$\bfi|_k$, $k\leq|\bfi|$, belongs to the tree;
\item if $\bfu\in\Up$, then $\bfu j\in\Up$ for $j = 1, \ldots,
Z_\bfu
$ and $\bfu j \notin\Up$ for $j > Z_\bfu$, where $Z_\bfu$ is the
number of children of $\bfu$.
\end{itemize}

Let $\Up_n$ be the $n$th generation of the tree, that is, the set of
nodes of length~$n$. (These are level $-n$ crossings, of size
$2^{-n}$.) Define\vspace*{1pt} $\Up_\bfi= \{\bfj\in\Up |
|\bfj|\geq|\bfi|$ and \mbox{$\bfj|_{|\bfi|}=\bfi\}$}. The boundary of
the tree is given by $\partial\Up= \{\bfi\in\N^\N |  \forall n\geq0,
\bfi|_n\in\Up\}$. Let $\psi(\bfi)$ be the position of node $\bfi$
within generation $|\bfi |$, so that crossing $\bfi$ is just
$C^{-|\bfi|}_{\psi(\bfi)}$. The nodes to the left and right of $\bfi$,
namely $C^{-|\bfi|}_{\psi (\bfi)-1}$ and $C^{-|\bfi|}_{\psi(\bfi)+1}$,
will\vspace*{1pt} be denoted $\bfi-$ and $\bfi+$. In general, when we
have quantities associated with crossings, we will use tree indexing
and level/position indexing interchangeably. So $Z_\bfi=
Z^{-|\bfi|}_{\psi(\bfi)}$, $T_\bfi= T^{-|\bfi|}_{\psi (\bfi )}$, etc.
At present\vspace*{1pt} our tree indexing only applies to crossings
contained within the first level 0 crossing; however, in Section
\ref{secextend} we will extend this notation to the whole tree.

Let $\al^n_k \in\{+, -\}$ be the orientation of $C^n_k$, $+$ for up
and $-$ for down, and let $A^n_k$ be the vector given by the
orientations of the subcrossings of $C^n_k$.
Let $D^n_k = T^n_k - T^n_{k-1}$ be the duration of $C^n_k$.
Clearly, to reconstruct the process we only need $\al^n_k$ and $D^n_k$
for all $n$ and $k$.
The $\al^n_k$ encode the spatial behavior of the process, and the
$D^n_k$ the temporal behavior.
Our definition of an EBP is concerned with the spatial component only.
%
\begin{definition}
A continuous process $X$ with $X(0) = 0$ is called an Embedded
Branching Process (EBP) process if for any fixed $n$, conditioned on
the crossing orientations $\al^n_k$, the random variables $A^n_k$ are
all mutually independent, and $A^n_k$ is conditionally independent of
all $A^m_j$ for $m > n$.
In addition we require that $\{A^n_k  |  \al^n_k = i\}$ are
identically distributed, for $i = +, -$.

That is, an EBP process is such that if we take any given crossing,
then count the orientations of its subcrossings at successively finer
scales, we get a (supercritical) two-type Galton--Watson process, where
the types correspond to the orientations.
\end{definition}

Subcrossing orientations have a particular structure.
A level $n$ up crossing is from $k2^n$ to $(k+1)2^n$, a down crossing
is from $k2^n$ to $(k-1)2^n$.
The level $n-1$ subcrossings that make up a level $n$ parent crossing
consist of \textit{excursions}\vadjust{\goodbreak} (up--down and down--up pairs) followed by a
\textit{direct crossing} (down--down or up--up pairs), whose direction
depends on the parent crossing: if the parent crossing is up, then the
subcrossings end up--up; otherwise, they end down--down.
Let $Z^n_k$ be the length of $A^n_k$, that is, the number of
subcrossings of $C^n_k$.
The number of up and down subcrossings will be written $Z^{n+}_k$ and
$Z^{n-}_k$, respectively.
Clearly, each of the $Z_k^n-2$ first entries of $A_k^n$ comes in pairs,
each pair being up--down or down--up. The last two components are either
the pair up--up or down--down, depending on~$\alpha_k^n$.
Thus, given $\al^n_k = +$, we must have $Z^{n+}_k = \half Z^n_k + 1$
and $Z^{n-}_k = \half Z^n_k - 1$, and conversely given $\al^n_k = -$.

Let $\calA$ be the space of possible orientations.
That is, $a \in\calA$ consists of some number of pairs, $+-$ or $-+$,
then a single pair $++$ or $--$.
Given an EBP process, for the offspring type distributions we write
$p^+_{A}(a) = \Pb(A^n_k = a  |  \al^n_k = +)$ and $p^-_{A}(a) = \Pb
(A^n_k = a  |  \al^n_k = -)$, for $a \in\calA$.
Let $\mu^+ = \E(Z^n_k | \al^n_k=+)$, $\mu^- = \E(Z^n_k  |  \al
^n_k=-)$ and $\mu= \half(\mu^+ + \mu^-)$, then the mean offspring
matrix is given by
\[
M:= \E\pmatrix{
(Z^{n+}_k|\al^n_k = +) & (Z^{n-}_k|\al^n_k = +) \cr
(Z^{n+}_k|\al^n_k = -) & (Z^{n-}_k|\al^n_k = -)}
=
\pmatrix{
\half\mu^+ + 1 & \half\mu^+ - 1 \cr
\half\mu^- - 1 & \half\mu^- + 1}.
\]

To proceed we need to make some assumptions about $p^\pm_{A}$.
%
\begin{assp}\label{AssGW}
$\mu^+, \mu^- > 2$ and $\E(Z^{ni}_k \log Z^{ni}_k  |  \al^n_k =
j) <
\infty$ for \mbox{$i, j = \pm$}.
\end{assp}

The first of\vspace*{1pt} these assumptions ensures that $M$ is
strictly positive with dominant eigenvalue $\mu> 2$, and corresponding
left eigenvector $(\half, \half)$. The corresponding\vspace*{1pt} right
eigenvector is $((\mu^+-2)/(\mu-2), (\mu ^--2)/(\mu-2))^T$. The second
assumption is the usual condition for the normed limit of a
supercritical Galton--Watson process to be nontrivial.
%
\begin{theorem}\label{CEBP}
For any offspring orientation distributions $p^\pm_{A}$ satisfying
Assumption~\ref{AssGW}, there exists a corresponding continuous EBP
process $X$ defined on~$\R_+$.
\end{theorem}
\begin{pf}
A version of this result can be found as Theorem 1 in~\cite{J04}, for
particular orientation distributions.

\begin{figure}[b]

\includegraphics{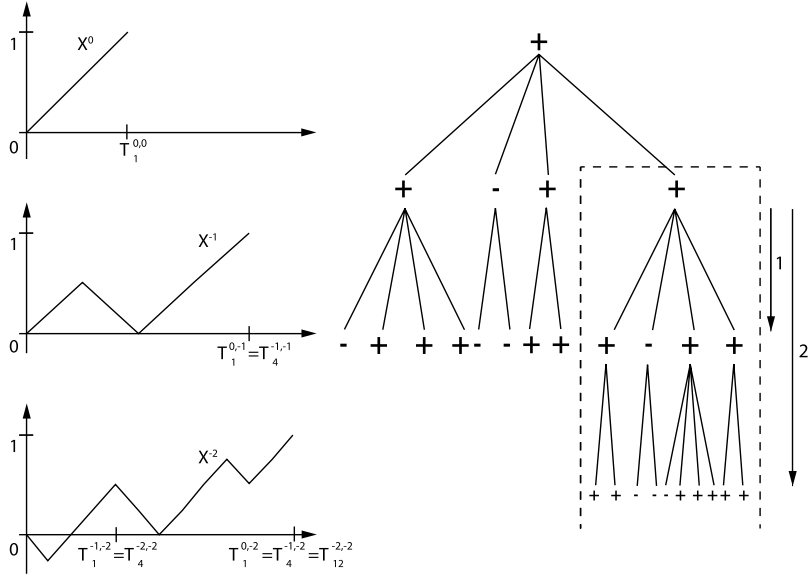}
 \caption{Construction of $X^0$, $X^{-1}$ and $X^{-2}$,
and the associated crossing tree (see the proof of Theorem
\protect\ref{CEBP}). The subtree rooted at crossing $C^{-1}_4$ (the 4th
crossing of size $1/2$) corresponds to the tree inside the dashed box.
In the notation of the proof of Theorem \protect\ref{CEBP}, for this
subtree we have $m=1$. If we go down one level in the subtree,
corresponding to level $n=2$ of the original tree, then
$S^{+}_{-1,4}(1) = 3$ and $S^{-}_{-1,4}(1) = 1$ count the
number\vspace*{-1pt} of up and down crossings at level 1 of the
subtree. Similarly, $S^{+}_{-1,4}(2) = 7$ and $S^{-}_{-1,4}(2) = 3$
count the number of up and down crossings at level 2 of the subtree,
and so on. This figure also illustrates other notation used in the
proof of Theorem~\protect\ref{CEBP}. For example, one has
$T^{-1,-2}_1=T^{-2,-2}_4$, since for $X^{-2}$, the the first crossing
time of size $2^{-1}$ corresponds to the fourth
crossing time of size $2^{-2}$.}%
\label{Spmfig}
\end{figure}

\textsc{Step 1.}
We initially construct a single crossing from 0 to 1, with support $[0, T^0_1]$.
In step 2 we will extend the range to $\R$ and the support to
$[0,\infty)$.
$X$~is obtained as the limit as $n \to+\infty$, of a sequence of
random walks $X^{-n}$ with steps of size $2^{-n}$ and duration $\mu^{-n}$.
Put $X^0(0)=0$ and $X^0(1)=1$, so that the coarsest scale is $n=0$.
Given $X^{-n}$ we construct $X^{-(n+1)}$ by replacing the $k$th step of
$X^{-n}$ by a sequence of $Z_k^{-n}$ steps of size $2^{-(n+1)}$ and
duration $\mu^{-(n+1)}$.
If $\al^{-n}_k = i$, then the orientations $A^{-n}_k$ of the
subcrossings are distributed according to $p^i_{A}$.
For a given $n$ the $A^{-n}_k$ are all mutually independent, and, given
$\al^{-n}_k$, $A^{-n}_k$ is conditionally independent of all
$A^{-m}_j$, for $-m > -n$.\vadjust{\goodbreak}

Denote the (random) time that $X^{-n}$ hits 1 by
\[
T^{0,-n}_1 = \inf\{t | X^{-n}(t)=1 \}.
\]
We define $X^{-n}(t)$ for all $t \in[0, T^{0,-n}_1]$ by linear
interpolation, and set $X^{-n}(t) = 1$ for all $t > T^{0,-n}_1$.
The interpolated $X^{-n}$ have continuous sample paths, and we will
show that they converge uniformly on any finite interval, from which
the continuity of the limit process follows.
For any $m \leq n$, let $T_0^{-m,-n}=0$ and
\[
T_{k+1}^{-m,-n}= \inf\{t> T_k^{-m,-n} | X^{-n}(t)\in2^{-m}\Z
,X^{-n}(t)\neq X^{-n}(T_k^{-m,-n})\}.
\]
If $X^{-n}(T_k^{-m,-n})=1$, then set $T_{k+1}^{-m,-n}= \infty$.
By construction $X^{-n}(T^{-m,-n}_k) = X^{-m}(\mu^{-m}k)$, for all $k$
and $m \leq n$.
The duration of the $k$th level $-m$ crossing of $X^{-n}$ is
$D_k^{-m,-n}= T_k^{-m,-n}- T_{k-1}^{-m,-n}$.

A realization of $X^0$, $X^{-1}$ and $X^{-2}$ is given in
Figure~\ref{Spmfig}, with the associated crossing tree.

We use a branching process result to establish that the crossing
durations converge.
When we defined the crossing tree (see Figure~\ref{crossingtree}) we
started with a sample path and then defined generations of crossings:
taking the first crossing of size 1 as the root (level or generation
0), its subcrossings of size $1/2$ form the second generation (or
level), its subcrossings of size $1/4$ form the third generation, and
so on.
Each crossing can be up or down, so our tree has two types of nodes.
Here we are reversing that process.
That is, we are growing a tree using a two-type Galton--Watson process,
and from the tree, constructing a sample path.
The offspring distributions for our tree are just $p^\pm_A$.
Given the tree at generation $n$, we get an approximate sample path by
taking a sequence of up and down steps of size $2^{-n}$ and duration
$\mu^{-n}$, with directions taken from the node types of the tree.
We need to show that the sequence of sample paths, obtained as $n\to
\infty$, converges.

Consider the subtree descending from crossing $C^{-m}_k$.
Let $S^{+}_{-m,k}(n-m)$ and $S^{-}_{-m,k}(n-m)$ be the number of up and
down crossings of size $2^{-n}$ which are descended from the $k$th
crossing of size $2^{-m}$; then $\{(S^{+}_{-m,k}(n-m),
S^{-}_{-m,k}(n-m))\}_{n=m}^\infty$ is a two-type Galton--Watson process.
From Athreya and Ney~\cite{AN72}, Section~V.6, Theorems 1 and 2, we
have that as $n \to\infty$, $\mu^{m-n}(S^{+}_{-m,k}(n-m),
S^{-}_{-m,k}(n-m))$ converges almost surely and in mean to $(\half,
\half) W^{-m}_k$, where $W^{-m}_k$ is strictly positive, continuous and
$\E(W^{-m}_k  |  \al^{-m}_k = \pm) = (\mu^\pm- 2)/(\mu- 2)$.
Moreover, the distribution of $W^{-m}_k$ depends only on $\al^{-m}_k$,
and for any fixed $m$ the $W^{-m}_k$ are all independent.
Finally, since $S^{+}_{-m,k}(n-m) + S^{-}_{-m,k}(n-m) = \mu
^{n}D^{-m,-n}_k$, we have
\[
D_k^{-m,-n}\rightarrow\mu^{-m} W_k^{-m} \qquad\mbox{a.s. as
}n\rightarrow
\infty.
\]
Accordingly, let $T_k^{-m} = \sum_{j=1}^k \mu^{-m} W_j^{-m} = \lim
_{n\rightarrow\infty} T_k^{-m,-n}$.\vspace*{1pt}

Take any $\varepsilon>0$, $\delta>0$ and $T>0$.
To establish the a.s. convergence of the processes $X^{-n}$, uniformly
on compact intervals, we show that we can find a $u$ so that with
probability $1-\varepsilon$,
%
\begin{equation}\label{eqconti}
|X^{-r}(t)- X^{-s}(t)|\leq\delta\qquad\mbox{for all } r,s \geq u
\mbox{ and } t\in[0,T].
\end{equation}
Given $t\in[0,T]$, let $k=k(n,t)$ be such that
\[
T_{k-1}^{-n}\leq t < T_k^{-n}.
\]
For any $r,s \geq n$, the triangle inequality yields
%
\begin{eqnarray}\label{eqineqtri}
|X^{-r}(t)- X^{-s}(t)|&\leq&|X^{-r}(t)-X^{-r}(T_k^{-n,-r})|\nonumber\\[-8pt]\\[-8pt]
&&{}+ |X^{-s}(T_k^{-n,-s})-X^{-s}(t)|,\nonumber
\end{eqnarray}
since $X^{-r}(T_k^{-n,-r})=X^{-s}(T_k^{-n,-s})= X^{-n}(k\mu^{-n})$.

For any $u \geq n$ let $j=j(n,u)$ be the smallest $j$ such that
$T_j^{-n,-u}> T$.
As $u\rightarrow+\infty$, $j(n,u)\rightarrow j(n) < \infty$ a.s., so
for any $n$ we can choose $\varepsilon_0$ such that
\[
\Pb\Bigl(\min _{i\leq j(n)}\mu^{-n} W_i^{-n} \geq\varepsilon_0
\Bigr) \geq1-\varepsilon,
\]
and $u$ such that for all $q \geq u$,
\[
\Pb\Bigl(\max _{i\leq j(n)}|T_i^{-n,-q}-T_i^{-n}| <\varepsilon_0
\Bigr) \geq1-\varepsilon,
\]
which yields
\[
\Pb\Bigl(\max _{i\leq j(n)}|T_i^{-n,-q}-T_i^{-n}|
< \min _{i\leq j(n)}\mu^{-n}W_i^{-n} \Bigr) \geq1 -
\varepsilon.
\]
Thus, given $n$ we can find $u$ such that for all $q \geq u$, with
probability at least $1-\varepsilon$,
\[
T^{-n,-q}_{k-2}<t<T_{k+1}^{-n,-q} \qquad\mbox{for all }t\in[0,T].
\]

Now, since $X^{-q}(T_{k-2}^{-n,-q})= X^{-n}((k-2)\mu^{-n})$,
$X^{-q}(T_{k+1}^{-n,-q})= X^{-n}((k+1)\mu^{-n})$, and in three steps
$X^{-n}$ can move at most distance $3\cdot2^{-n}$, we have
\[
|X^{-q}(t)-X^{-q}(T_k^{-n,-q})|\leq3\cdot2^{-n}.
\]
Choosing $n$ large enough that $6\cdot2^{-n}\leq\delta$, we see that
(\ref{eqconti}) follows from (\ref{eqineqtri}).
Sending $\delta$ and $\varepsilon$ to $0$ shows that $X^{-n}$ converges to
some continuous limit process $X$ uniformly on all closed intervals
$[0, T]$, with probability $1$.
By construction, the duration of crossing $C^{-n}_k$ is $\mu^{-n}W^{-n}_k$.

\textsc{Step 2.}
Clearly the construction above can be used to generate any crossing
from 0 to $\pm2^n$.
Thus, to extend our construction from a single crossing to a process
$X(t)$ defined for all $t \in\R_+$, we proceed by constructing a
nested sequence of processes $\{X^{(n)}\}_{n=0}^\infty$, such that
$X^{(n)}$ is a crossing from 0 to $\pm2^n$, and the first level $n$
crossing of $X^{(n+1)}$ is precisely $X^{(n)}$.
To make this work, we just need to specify $\Pb(X^{(n)}(T^n_1) = 2^n)$
in a consistent manner.

Consider the orientation of the first crossing from 0 to $\pm2^n$ for
an EBP process.
Let $u = \Pb(\al^n_1 = +  |  \al^{n+1}_1 = +)$ and $v = \Pb(\al
^n_1 =
+  |  \al^{n+1}_1 = -)$; then $u$ and $v$ are determined by $p^\pm
_{A}$, and
%
\begin{equation}\label{firstcrossingeqn}
a_n:= \Pb(\al^n_1 = +) = u a_{n+1} + v (1 - a_{n+1}) = v + (u-v) a_{n+1}.
\end{equation}
For $(u, v) \in[0, 1]^2 \setminus\{(1,0)\}$, we see that equation
(\ref{firstcrossingeqn}) has fixed point $a = v/(1-u+v) \in[0, 1]$.
Moreover, the only doubly infinite sequence $\{ a_n \}_{n = -\infty
}^\infty$ which satisfies (\ref{firstcrossingeqn}) and remains in
$[0, 1]$ is given by $a_n = a$ for all $n$.
Given this, it follows that $a_n = a$, and thus from Bayes's theorem
that $\Pb(\al^{n+1}=+  |  \al^n_1=+) = u$ and $\Pb(\al^{n+1}=+
|
\al^n_1=-) = v$.
If $(u, v) = (1, 0)$, then any $a \in[0, 1]$ is possible, but
everything else goes through as before.
In this case the $\al^n_1$ are all the same, but may be of either type.

Construct $X^{(0)}$ as a crossing from 0 to 1 with probability $a$ [the
fixed point of~(\ref{firstcrossingeqn})], otherwise as a
crossing from 0 to $-1$.
Then, given $X^{(n)}$, construct $X^{(n+1)}$ as follows: first, put
$\al
^{n+1}_1 = +$ with probability $u$ if $\al^0_1=+$, with probability $v$
otherwise; second, generate $A^{n+1}_1$ conditional on $\al^{n+1}_1$
and $\al^n_1$; third, use $X^{(n)}$ as the first level $n$ crossing of
$X^{(n+1)}$; finally construct the remaining level $n$ crossings
conditional on $\al^{n}_2, \al^{n}_3, \ldots, \al^{n}_{Z^{n+1}_1}$.
Write $X$ for the limit of the $X^{(n)}$.
To complete our construction we just need to check that the process $X$
does not escape to $\pm\infty$ in finite time.
By construction, we have $T^n_1 = \inf\{ t  |  X(t) = \pm2^n\} =
\mu
^n W^n_1$, where $W^n_1$ is strictly positive, continuous, and has a
distribution depending only on the orientation $\al^n_1$.
Thus for any $T > 0$, $\Pb(T^n_1 < T) \to0$ as $n\to\infty$.
\end{pf}
%
\begin{theorem}\label{CEBP2}
Let $X$ be the EBP constructed in Theorem~\ref{CEBP}; then, for
each~$n$, conditioned on the crossing orientations $\al^n_k$, the
crossing durations $D^n_k$ are all mutually independent, and $D^n_k$ is
conditionally independent of all $A^m_j$ for $m > n$. Also, $\E(D^n_k
|  \al^n_k = \pm) = \mu^{n} (\mu^\pm- 2)/(\mu- 2)$, and the
distribution of $\mu^{-n}D^n_k$ depends only on $\al^n_k$. Moreover, up
to finite-dimensional distributions, $X$~is the unique such EBP with
offspring orientation distributions $p^\pm_{A}$. That is,\vspace*{1pt} for any other
EBP process $Y$ with offspring orientation distributions $p^\pm_A$ and
crossing durations as above, we have $(X(t_1), \ldots, X(t_k))
\eqdist(Y(t_1), \ldots, Y(t_k))$ for any $0 \leq t_1 < t_2 < \cdots<
t_k$.

Accordingly, we call $X$ the Canonical EBP (CEBP) process with these
offspring distributions.

We also observe that $X$ is discrete scale-invariant: let $H=\log
2/\log
\mu$; then for all $c\in\{\mu^n, n\in\Z\}$,
%
\begin{equation}\label{eqholder}
X(t) \stackrel{\mathit{fdd}}{=} c^{-H} X(ct),
\end{equation}
where $\stackrel{\mathit{fdd}}{=}$ denotes equality for finite-dimensional
distributions.
$H = \log\mu/\log2$ is known as the Hurst index.
\end{theorem}
\begin{pf}
We retain the notation of Theorem~\ref{CEBP}.

For the process $X$, the dependence structure of the crossing durations
is clear from the construction.

To show\vspace*{1pt} uniqueness, let $Y$ be some other EBP process with offspring
orientation distributions $p^\pm_{A}$, and crossing durations
satisfying the conditions of the theorem statement.
We will make use of the same notation for the crossing times,
durations, orientations, etc. of $Y$ as for $X$, and rely on the
context to distinguish them.

For an EBP, the finite joint distributions of the orientations $A^n_k$
are determined completely by $p^\pm_{A}$, and thus are identical for
$X$ and $Y$.\vadjust{\goodbreak}
For the crossing durations of $Y$, note that for any $m \leq n$ and
$k$, we have
%
\begin{equation}\label{poinc}
\mu^m D^{-m}_k = \mu^{m-n} \sum_{j=\zeta(-m,-n,k)+1}^{\zeta(-m,-n,k+1)}
\mu^{n} D_{j}^{-n},
\end{equation}
where $\zeta(-m,-n,k)$ is such that $\zeta(-m,-n,k)+1$ is the index of
the first level $-n$ subcrossing of $C^{-m}_k$.
Thus by the strong law of large numbers, sending $n\to\infty$,
\begin{eqnarray*}
\mu^m D^{-m}_k
&=& \mu^{m-n}S^{+}_{-m,k}(n-m)
\Biggl[ \frac{\mu^n}{S^{+}_{-m,k}(n-m)}
\mathop{\sum_{j=\zeta(-m,-n,k)+1}}_{\al^{-n}_j=+}^{\zeta(-m,-n,k+1)}
D^{-n}_j \Biggr] \\[-2pt]
&&{}+ \mu^{m-n}S^{-}_{-m,k}(n-m)
\Biggl[ \frac{\mu^n}{S^{-}_{-m,k}(n-m)}
\mathop{\sum_{j=\zeta(-m,-n,k)+1 }}_{
\al ^{-n}_j=-}
^{\zeta(-m,-n,k+1)}
D^{-n}_j \Biggr] \\[-2pt]
&\convergeprob& \half W^{-m}_k \mu^n \E(D^{-n}_j  |  \al
^{-n}_j=+) +
\half W^{-m}_k \mu^n \E(D^{-n}_j  |  \al^{-n}_j=-) \\[-2pt]
&=& W^{-m}_k,
\end{eqnarray*}
where the distribution of $W^{-m}_k$ is completely determined by $p^\pm
_{A}$, and thus is the same for $X$ and $Y$.

Once we have the crossing orientations and the assumed dependence
structure of the crossing durations, the crossing distributions (for up
and down types) determine the joint distributions of the crossing times
$\{ T^n _k \}$.
Thus, for any $n$ and $k$, $\{ X(T^n_i) \}_{i=0}^k$ and $\{ Y(T^n_i) \}
_{i=0}^k$ are identically distributed.
Since any $t$ can be bracketed by a sequence of hitting times, $X$ and
$Y$ are identical up to finite-dimensional distributions.

That $X$ is discrete scale-invariant is a direct consequence of its
construction, since simultaneously scaling the state space by $2^k$ and
time space by $\mu^k$ does not change the distribution of $X$.
\end{pf}
%
\begin{remark}
From~\cite{H92} it is clear that Brownian motion is an example of a
CEBP process, where the offspring of any crossing consist of a
geometric ($1/2$) number of excursions, each up--down or down--up with
equal probability, followed by either an up--up or down--down direct crossing.
That is,\looseness=-1
\[
p^+_{A}(\cdots{+} {+}) = p^-_{A}(\cdots{-} {-}) = 2^{-(z+1)},
\]\looseness=0
where $\cdots$ represents a combination of $z$ pairs, each either $+ -$
or $- +$.
It follows that $\Pb(Z^n_k = 2x) = 2^{-x}$, independently of $\al^n_k$.
\end{remark}

\section{From CEBP to MEBP}
\label{CEBPMEBP}

In this section we construct Multifractal Embedded Branching processes
(MEBP processes) as time changed CEBP processes.\vadjust{\goodbreak}

Consider initially a single crossing of a CEBP $X$, from 0 to $\pm1$.
We constructed $X$ as the limit of a sequence of processes $X^{-n}$,
which take steps of size $2^{-n}$ and duration $\mu^{-n}$.
The crossing tree gives the number of subcrossings of each crossing.
If we add a weight of $1/\mu$ to each branch of the tree, then
truncating the tree at level $-n$, the product of the weights down any
line of descent is $\mu^{-n}$, which is the duration of any single
crossing by $X^{-n}$.
We generalize this by allowing the weights to be random, then defining
the duration of a crossing to be the product of the random weights down
the line of descent of the crossing.
The resulting process, $Y^{-n}$ say, can be viewed as a time-change of
$X^{-n}$, where the time-change is obtained from a multiplicative
cascade defined on a (two-type) Galton--Watson tree.

As for CEBP, we will initially construct a single level 0 crossing of
an MEBP, then extend the construction to $\R_+$.
We will retain the notation of Section~\ref{secEBP}, but note that we
will prefer the tree indexing scheme to the level/position indexing
scheme in what follows.
In particular, the number of level $-n$ up and down subcrossings of
node $\bfi$ in level $-m$ are denoted $S^+_\bfi(n-m)$ and $S^-_\bfi
(n-m)$, and, under Assumption~\ref{AssGW}, the almost sure limit and
mean limit of $\mu^{m-n}(S^+_\bfi(n-m), S^-_\bfi(n-m))$ is $(\half,
\half)W_\bfi$.
The duration of crossing $\bfi$ of the CEBP process $X$ is then $\mu
^{-m}W_{\bfi}$.

We assign weight $R_\bfi(j)$ to the branch $(\bfi, \bfi j)$.
$R_\bfi:= (R_\bfi(1), \ldots, R_\bfi(Z_\bfi))$ may depend on
$A_\bfi$,
but conditioned on $\al_\bfi$ must be independent of other nodes
$\bfj$
that are not descendants of $\bfi$.
For $r \in\R_+^{|a|}$, write $F^\pm_{R|a}(r) = \Pb(R_\varnothing(1)
\leq
r(1), \ldots, R_\varnothing(z) \leq r(z)  |  \al_\varnothing=\pm,
A_\varnothing=a, |a|=z)$ for the joint distribution of $R_\varnothing$,
conditioned on the crossing orientations $a$.
The weight attributed to node $\bfi$ is
\[
\rho_{\bfi} = \prod _{k=0}^{|\bfi|-1} R_{\bfi|_{k}}(\bfi[k+1]).
\]
That is, $\rho_\bfi$ is the product of all weights on the line of
descent from the root down to node $\bfi$.
We use the weights to define a measure, $\nu$, on the boundary of the
crossing tree.
The measure $\nu$ on $\partial\Up$ is then mapped to a measure
$\zeta$
on $\R$, with which we define a chronometer $\calM$ (a nondecreasing
process) by\vspace*{1pt} $\calM(t) = \zeta([0, t])$.
The MEBP process is then given by $Y = X \circ\calM^{-1}$, where $X$
is the CEBP.
The crossing trees of $X$ and $Y$ have the same spatial structure, but
have different crossing durations.
In Figure~\ref{EBPMEBP} we plot a realization of an MEBP process and
its associated CEBP.

The literature on multiplicative cascades is rather extensive.
For the existence of limit random measures and the study of the
properties of certain martingales defined on $m$-ary trees, one can
refer, for instance, to the works of Kahane and Peyri\`ere~\cite{KP76},
Barral~\cite{B99}, Liu and Rouault~\cite{LR00} and Peyri\`ere~\cite{P00}.
For results on random cascades defined on Galton--Watson trees, see,
for example, Liu~\cite{L99,L00}, and Burd and Waymire~\cite{BW00}.

To obtain the time-change process explicitly, the random measure
defined on the boundary of the tree is mapped to $\R_+$ then
integrated. Note that this mapping, given explicitly in Section
\ref{cascade}, differs from random partitions previously considered in
the literature. The usual approach is to use a ``splitting measure'' to
map the boundary of the tree to $[0, 1]$, then use the density of the
cascade measure with respect to the splitting measure; see, for
example,~\cite{P77,P79,R03}. Our approach can be thought of as using
the ``branching measure'' instead of a splitting measure. A~splitting
measure is constructed by splitting the mass associated with a given
node between its offspring, with no mass lost or gained. The branching
measure allocates mass according to the number of offspring, and is
only conserved in mean. We have taken the terminology of splitting and
branching measures from~\cite{L00}, Example 1.3. A multifractal study
of the measure we construct on $\R_+$ is given in a forth-coming paper
\cite{DHJ}.

\subsection{\texorpdfstring{The measure $\nu$}{The measure nu}}\label{measurenu}

To construct $\nu$, we use a well-known correspondence between
branching random walks and random cascades, in which the offspring of
individual $\bfi$ have types given by $A_\bfi$ and displacements
(relative to $\bfi$) given by $-{\log R_\bfi}$.
For background on multitype branching random walks, we refer the reader
to Kyprianou and Sani~\cite{KS01} and Biggins and Sani~\cite{BS05}.

Suppose $|\bfi|=m$ and $n \geq m$.
Define
\[
\calW^\pm_\bfi(n-m, \tha) = \sum_{\bfj\in\Up_n\cap\Up_\bfi,
\al_\bfj=\pm
} (\rho_\bfj/\rho_\bfi)^\tha
\]
and for $i,j = \pm$,
\begin{eqnarray*}
m_{i,j}(\tha) &=& \E\bigl(\calW^j_\varnothing(1, \tha)  |  \al
_\varnothing=i\bigr)
\\
&=& \E\biggl(\sum_{1\leq k\leq Z_\varnothing, \al_k=j}
R_\varnothing
(k)^\tha \Big|  \al_\varnothing=i \biggr).
\end{eqnarray*}
Let $M(\tha) = (m_{i,j}(\tha) )_{i,j=\pm}$, and write
$m^n_{i,j}(\tha)$
for the $(i,j)$ entry of the $n$th power $M^n(\tha)$. Then it is
straight forward to check that $\E( \calW^j_\bfi(n-m, \tha)  |
\al
_\bfi=i) = m^{n-m}_{i,j}(\tha)$.
If we take constant weights equal to $1/\mu$, then $\calW^\pm_\bfi(n-m,
1) = \mu^{m-n} S^\pm_\bfi$, in the notation of Theorem~\ref{CEBP}.

Let $\mu(\tha)$ be the largest eigenvalue of $M(\tha)$.
We make the following assumptions about $R_\varnothing$.
%
\begin{assp}\label{Ass1}
We suppose that $0 < R_\varnothing< \infty$ a.s., $M(\tha) < \infty$ in
an open neighborhood of 1, $\mu(1) = 1$ and $\mu'(1) < 0$.

In the case where the distribution of $R_\varnothing$ (and thus
$Z_\varnothing$) does not depend on $\al_\varnothing$, we assume in
addition that
\[
\E\Biggl( \sum_{j=1}^{Z_\varnothing} R_\varnothing(j) \log\sum
_{j=1}^{Z_\varnothing} R_\varnothing(j) \Biggr) < \infty.
\]

In the case where there is dependence on the crossing orientation
(type), we suppose that for some $\delta> 1$, $\mu(\delta) < 1$ and
\[
\E\Biggl(  \Biggl( \sum_{j=1}^{Z_\varnothing} R_\varnothing(j)
\Biggr)^\delta \bigg|  \al_\varnothing=i \Biggr) < \infty\qquad\mbox{for }
i =
\pm.
\]
\end{assp}

Note that if the weights are finite and strictly positive, then $M$ and
$\mu$ from the previous section are just $M(0)$ and $\mu(0)$, and from
Assumption~\ref{AssGW} we get $0 < M(\tha)$ for all $\tha\geq0$.
In the case where $R_\varnothing$ does not depend on $\al_\varnothing$, the
BRW simplifies to a single-type process, and the condition $\mu(1) = 1$
simplifies to $\E( \sum_{j=1}^{Z_\varnothing} R_\varnothing(j)
)
= 1$, which we recognize as a conservation of mass condition.

Left and right eigenvectors corresponding to $\mu(1)$ will be denoted
$\bfu= (u^+, u^-)$ and $\bfv= (v^+, v^-)^T$, normed so that $\bfu(1,
1)^T = 1$ and $\bfu\bfv= 1$.
The following lemma is a direct consequence of Biggins and Kyprianou
\cite{BK04}, Theorem 7.1, and Biggins and Sani~\cite{BS05}, Theorem 4.
%
\begin{lem}\label{BRWlem}
Under Assumptions~\ref{AssGW} and~\ref{Ass1}, $(\calW^+_\bfi(n-m, 1),
\calW^-_\bfi(n-m, 1))$ converges almost surely to $\bfu\calW_\bfi$, for
some random variable $\calW_\bfi$ such that the distribution of
$\calW
_\bfi$ depends only on $\al_\bfi$, and $\E( \calW_\bfi |  \al
_\bfi=
i) = v^i$.
Moreover, for each $n$, conditioned on the crossing orientations $\al
_\bfi$, $\bfi\in\Up_n$, the $\calW_\bfi$ are mutually independent,
and $\calW_\bfi$ is conditionally independent of $(A_\bfj, R_\bfj)$ for
$|\bfj| < |\bfi|$.
For all nodes~$\bfi$,
%
\begin{equation}\label{eqE}
\calW_\bfi= \sum _{j=1}^{Z_{\bfi}} R_\bfi(j) \calW_{\bfi j}.
\end{equation}
\end{lem}

Note that in the case where $R_\varnothing$ does not depend on $\al
_\varnothing$, the right eigenvector $\bfv= (1, 1)^T$.

We can now define the measure $\nu$ on $\dUp$.
Recall $\Up_\bfi= \{\bfj\in\Up |  |\bfj|\geq|\bfi|$
and
\mbox{$\bfj|_{|\bfi|}=\bfi\}$}, so $\dUp_{\bfi}$ contains all the nodes on the
boundary of the tree which have $\bfi$ as an ancestor.
We define $\nu(\dUp_\bfi) = \rho_\bfi\calW_\bfi$.
By Carath\'eodory's extension theorem, we can uniquely extend $\nu$ to
the sigma algebra generated by these cylinder sets.

\subsection{\texorpdfstring{The measure $\zeta$ and time change $\calM$}
{The measure zeta and time change $\calM$}}\label{cascade}

The measure $\zeta$ is a mapping of $\nu$ from $\dUp$ to $[0,
W_\varnothing] \subset\R$.
By analogy with $m$-ary cascades, we call $\zeta$ a Galton--Watson
cascade measure on $[0, W_\varnothing]$.

As above, let $T^{-n}_k$ denote the $k$th level $-n$ passage time of
the CEBP process~$X$, and put
\[
\zeta((T^{-n}_{k-1}, T^{-n}_k]):= \nu(\dUp^{-n}_k)= \rho^{-n}_k
\calW^{-n}_k.
\]
Putting $\zeta(\{0\}) = 0$, this gives us $\zeta([0, T^{-n}_k])$ for
all $n, k \geq0$.
For arbitrary $t \in(0, W_\varnothing]$, let $\bfi\in\dUp$ be such
that $t \in(T^{-n}_{\psi(\bfi|_n)-1}, T^{-n}_{\psi(\bfi|_n)}]$ for all
$n \geq0$.
Noting that $T^{-n}_{\psi(\bfi|_n)} = T_{\bfi|_n}$ is a nonincreasing
sequence, we define $\zeta([0, t]) = \lim_{n\to\infty}\zeta([0,
T^{-n}_{\psi(\bfi|_n)}])$.\vspace*{2pt}

We can now define $\calM(t) = \zeta([0,t])$, and define the MEBP
process $Y$ (on $[0, \calW_\varnothing]$) as
\[
Y = X \circ{\cal M}^{-1}.
\]
Here we take $\calM^{-1}(t) = \inf\{s\dvtx\calM(s) \geq t\}$, so that
it is well defined, even if $\calM$ has jumps or flat spots.

Put $\calT^{-n}_k = \calM(T_{k}^{-n}) = \sum_{j=1}^k \rho^{-n}_j
\calW
^{-n}_j$. Then $Y({\cal T}_{k}^{-n}) = X(T_{k}^{-n})$, so $\calT
^{-n}_k$ is the $k$th level $-n$ crossing time for $Y$, and $\calD
^{-n}_k = \rho^{-n}_k \calW^{-n}_k$ the $k$th level $-n$ crossing duration.
Note that if we take constant weights equal to $1/\mu$, then $\calT
^{-n}_k = T_{k}^{-n}$ and $Y=X$.
%
\begin{lem}
Under Assumptions~\ref{AssGW} and~\ref{Ass1}, $\calM$ and $\calM^{-1}$
are continuous.
That is, $\calM$ has neither jumps nor flat spots.
\end{lem}
\begin{pf}
To show that $\calM$ has no flat spots, it is enough to show that:
\begin{longlist}[(b)]
\item[(a)]
\[
\max_k \mu^{-n} W^{-n}_k \convergeprob0 \qquad\mbox{as } n\to\infty,
\]
\item[(b)]
\[
\calW^{-n}_k > 0 \qquad\mbox{a.s. for each } n, k \geq0.
\]
\end{longlist}

Property (a) follows directly from Theorem 1 in~\cite{O80}, noting that
under Assumption~\ref{AssGW} $\E( W_\bfi |  \al_\bfi=\pm) <
\infty$,
so that $\int_0^y x \,dF^\pm(x)$ is slowly varying, where $F^\pm(x) =
\Pb
( W_\bfi\leq x  |  \al_\bfi=\pm)$.
This is equivalent to saying that the measure $\bar{\nu}$, defined on
$\dUp$ by $\bar{\nu}(\dUp_\bfi) = \mu^{-|\bfi|}W_\bfi$, has no atoms.

To show (b), let $q_\pm= \Pb(\calW_\varnothing=0  |  \al_\varnothing
=\pm
)$, then note that since the weights $R_\varnothing> 0$, we have, from
(\ref{eqE}), that
%
\begin{equation}\label{eqE2}
q_i = f_{i}(q_+, q_-),
\end{equation}
where $f_{i}$ is the joint probability generating function of $Z^\pm
_\varnothing$ given $\al_\varnothing=i$.
[Note that $\bar{q}_\pm= \Pb(W_\varnothing=0  |  \al_\varnothing
=\pm)$
satisfy the same equations.]
Since $(Z^i_\varnothing |  \al_\varnothing=i) \geq2$ and $\Pb(Z^\pm
_\varnothing= 2  |  \al_\varnothing=i) < 1$, we have for $(q_+, q_-)
\in
[0, 1]^2 \setminus\{ (0,0), (1,1) \}$, $f_i(q_+, q_-) < q_i$.
Thus the only solutions to (\ref{eqE2}) are $(0, 0)$ and $(1, 1)$, and
as $\E( \calW_\varnothing |  \al_\varnothing=\pm) > 0$, we get
$q_\pm= 0$.

$\calM$ is continuous (has no jumps) if $\zeta$ has no atoms.
That is,
\begin{longlist}[(b*)]
\item[(a*)]
\[
\max_k \rho^{-n}_k \calW^{-n}_k \convergeprob0 \qquad\mbox{as } n\to
\infty,\vadjust{\goodbreak}
\]
\item[(b*)]
\[
W^{-n}_k > 0 \qquad\mbox{a.s. for each } n, k \geq0.
\]
\end{longlist}
We prove (b*) in exactly the same way as (b).

Property (a*) is equivalent to saying that $\nu$ has no atoms.
In the case where the distribution of $R_\varnothing$ does not depend on
$\al_\varnothing$, the BRW embedded in the crossing tree is effectively
single-type, and (a*) is given by Liu and Rouault~\cite{LR97}, Theorem 6.
In the case where the distribution of $R_\varnothing$ does depend on
$\al
_\varnothing$, the approach of~\cite{LR97} generalizes only as far as the
end of their Lemma 13, at which point we require, for some $\lambda< 1$,
%
\begin{equation}\label{rightmosteqn}
\E\nu(\{ \bfi\dvtx\rho_{\bfi|_n} \geq\lambda^n\}) \to0
\qquad\mbox{as }
n \to\infty.
\end{equation}
However, this can be shown using some recent results of Biggins \cite
{B10}, as we now demonstrate.

In the notation of~\cite{B10}, consider a BRW with offspring types $\{
\sigma_i \} \eqdist\{ A_\varnothing(i) \}$ and displacements $\{ z_i \}
\eqdist\{ \log(R_\varnothing(i)\gamma) \}$, for some $\gamma> 1$.
Put
\[
\bar{m}_{i,j}(\theta) = \E\biggl( \sum_{1 \leq k \leq
Z_\varnothing,
A_\varnothing(k)=j} R_\varnothing(k)^\theta\gamma^\theta \Big|  \al
_\varnothing=i \biggr)
\]
(this is $m_{i,j}$ in the notation of \cite
{B10}). Then the matrix $\Mbar(\theta) = (\bar{m}_{i,j}(\theta
))_{i,j=\pm
}$ has maximum ``Perron--Frobenius'' eigenvalue $\kappa(\theta) = \mu
(\theta) \gamma^\theta$.
From assumptions~\ref{AssGW} and~\ref{Ass1} it is clear that for some
$\theta> 0$, $\Mbar(\theta)$ is finite, irreducible and primitive.

Let $\calB^{(n)}_i$ be the rightmost particle of type $i$ in generation
$n$, that is,
\[
\calB^{(n)}_\pm= \max_{\bfi\in\Up_n, \al_\bfi=\pm} \log\rho
_\bfi+ n
\log\gamma.
\]
Then Proposition 5.6 of~\cite{B10} shows that
\[
\frac{\calB^{(n)}_\pm}n \convergeas\Gamma(\kappa^*),
\]
where $\kappa^*(a) = \sup_{\theta\geq0}\{\theta a - \kappa(\theta
)\}$
and $\Gamma(\kappa^*) = \sup\{a\dvtx\kappa^*(a) < 0\}$.

We have $\kappa(0) = \mu(0)$, $\kappa(1) = \gamma$, $\kappa'(1) =
\mu
'(1)\gamma+ 1$, and for $\gamma$ large enough, $\kappa(\theta) \to
\infty$ as $\theta\to\infty$, faster than linear.
$\Gamma(\kappa^*)$ corresponds to the slope of the line that passes
through the origin and is tangent to $\kappa^*$, from which it follows
that $\Gamma(\kappa^*) < \gamma$ provided that $\kappa'(1) \neq
\gamma
$, that is, provided $\mu'(1) \neq(\gamma- 1)/\gamma$.
But $\mu'(1) < 0$ and $\gamma> 1$ by assumption, so $\kappa'(1) \neq
\gamma$, and we get
\[
\max_{\bfi\in\Up_n, \al_\bfi=\pm} \frac{\log\rho_\bfi}n
\convergeas
\Gamma(\kappa^*) - \gamma< 0.
\]
Equation (\ref{rightmosteqn}) follows immediately, completing the
proof of our lemma.
\end{pf}

\subsection{Extending the construction to $\R_+$}\label{secextend}

We can extend $Y$ from $[0, \calW_\varnothing]$ to $\R_+$ in much the
same way we extended the CEBP $X$, by constructing a sequence of nested
processes $Y^{(n)}$, where $Y^{(n)}$ consists of a a single level $n$\vadjust{\goodbreak}
crossing from 0 to $\pm2^n$, and the first level $n$ crossing of
$Y^{(n+1)}$ is precisely $Y^{(n)}$.
As for the CEBP we need to specify $\Pb(Y^{(n)}(\calT^n_1) = 2^n)$ in a
consistent manner, but we also need to scale the first crossing.

Construct $Y^{(0)}$ as a crossing from 0 to 1 with probability $a$ [the
fixed point of (\ref{firstcrossingeqn})], otherwise as a
crossing from 0 to $-1$.
Then, given $Y^{(n)}$, construct $Y^{(n+1)}$ as follows: first, put
$\al
^{n+1}_1 = +$ with probability $u$ if $\al^n_1=+$ and probability $v$
otherwise; second, generate $(A^{n+1}_1, R^{n+1}_1)$ conditional on
$\al
^{n+1}_1$ and $\al^n_1$; third, scale the weights $R^{n+1}_1$ by
$1/R^{n+1}_1(1)$; fourth, use $Y^{(n)}$ as the first level $n$ crossing
of $Y^{(n+1)}$; finally, construct the remaining level $n$ crossings
conditional on\vspace*{-2pt} $\al^n_2, \al^n_3, \ldots, \al^n_{Z^{n+1}_1}$.
Write $Y$ for the limit of the $Y^{(n)}$.

When constructing $Y^{(n+1)}$ we take $Z^{n+1}_1$ independent
processes, each constructed like $Y^{(n)}$, then scale the first by $1
= R^{n+1}_1(1)/R^{n+1}_1(1)$, the second by
$R^{n+1}_1(2)/R^{n+1}_1(1)$, and so on, before stitching them together.
When constructing the second and subsequent level $n$ crossings of
$Y^{(n+1)}$, we proceed exactly as for the construction of $Y^{(0)}$,
except for a spatial scaling of $2^n$ and a temporal scaling of $\prod
_{k=1}^n 1/R^{k}_1(1)$, noting that the $R^{k}_1(1)$ are taken from the
first level $n$ crossing, and are thus independent of the second and
subsequent level $n$ crossings.
Thus with this construction, the process $Y^{(n)}(t)$ is distributed as
$2^nY^{(0)}(t \rho^{-n}_1)$, where $\rho^{-n}_1$ is the weight given to
the first level $-n$ crossing of $Y^{(0)}$ (a product of $n$ weights,
from level $-1$ to $-n$).

To complete our construction, we just need to check that the process
$Y$ does not escape to $\pm\infty$ in finite time.
To see this note that the \textit{second} level $n$ crossing of
$Y^{(n+1)}$ is distributed as
\[
\frac{R^{n+1}_1(2)}{\prod_{k=1}^{n+1} R^{k}_1(1)} \calW^n_2,
\]
where, conditioned on its orientation, $\calW^n_2$ is equal in
distribution to the level 0 crossing of $Y^{(0)}$, and is independent
of $R^{k}_1(1)$ for $k = 1, \ldots, n+1$ and of $R^{n+1}_1(2)$.
We have already seen that $\calW^n_2 > 0$ almost surely, and by
assumption, $R^{n+1}_1(2) > 0$, so it suffices to show that $\prod
_{k=1}^{n+1} R^{k}_1(1) \to0$ almost surely as $n\to\infty$.
Given the orientations $\al^k_1$, $k = 1, \ldots, n+1$, the weights
$R^k_1(1)$ are independent.
The sequence of orientations $\{ \al^k_1 \}_{k=1}^\infty$ form a
two-state ($+$ and $-$) Markov chain, with transition matrix
\[
\pmatrix{u & 1-u \cr v & 1-v}.
\]
Thus the product $R^1_1(1) R^2_1(1) \cdots$ can be written as a product
of independent random variables of the form
\[
C = \prod_{k=1}^U A_k \prod_{k=1}^V B_k,
\]
where $U \sim\operatorname{geom}(u)$, $V \sim\operatorname{geom}(1-v)$, $A_k \sim
(R^k_1(1)|\al^k_1=+)$, $B_k \sim(R^k_1(1)|\break \al^k_1(1)=-)$, and they are
all independent.
The product $\prod_{k=1}^{n+1} R^{k}_1(1)$ converges to zero if the sum
$\sum_{k=1}^{n+1} \log R^{k}_1(1)$ diverges to $-\infty$, which follows
almost surely from the strong law of large numbers, provided $\E\log C
= \frac1{1-u} \E\log A_1 + \frac1v \E\log B_1 < 0$ (assuming $u\neq
1$ and $v\neq0$).
That is, the process $Y$ is defined on $\R_+$ provided the following
assumption holds.
%
\begin{assp}\label{assInf}
If $u = \Pb(\al^n_1 = +  |  \al^{n+1}_1 = +) \neq1$ and $v = \Pb
(\al
^n_1 = \break +  |  \al^{n+1}_1 = -) \neq0$, then we suppose that
\[
\frac1{1-u} \E\bigl(\log R^n_1(1) | \al^n_1=+\bigr) + \frac1v \E\bigl(\log R^n_1(1)
| \al^n_1=-\bigr) < 0.
\]
If $u = 1$, then we require $\E(\log R^n_1(1) | \al^n_1=+) < 0$, and if
$v = 0$, then we require $\E(\log R^n_1(1) | \al^n_1=-) < 0$.
\end{assp}

To describe the crossing tree of the extended process $Y$, it is
convenient to extend the tree-indexing notation introduced earlier.
We do this by indexing nodes relative to a \textit{spine}, defined by the
first crossing at each level.
For any node in the tree, we can trace its ancestry back to the spine.
For any $n$ let $n\dvtx\varnothing$ be the node on level $n$ of the spine
and $\Up_{n\dvtx\varnothing}$ the tree descending from that node.
Nodes in the tree $\Up_{n\dvtx\varnothing}$ will be labeled $n\dvtx\bfi$, where
$\bfi$ is the node index relative to $n\dvtx\varnothing$.
Thus $n\dvtx\bfi$ is in level $n-|\bfi|$ of the crossing tree, and a
crossing previously labeled $\bfi$ is now labeled $0\dvtx\bfi$.
Note that this labeling is not unique, as $n\dvtx\bfi= (n+1)\dvtx1\bfi$.

Write $\rho_{n\dvtx\bfi}$ for the weight assigned to node $n\dvtx\bfi$, which
is given by
\[
\rho_{n\dvtx\bfi} = \cases{
\displaystyle \prod_{k=0}^{|\bfi|-1} R_{n\dvtx\bfi|_{k}}(\bfi[k+1]) \bigg/\prod_{k=0}^{(n
\wedge|\bfi|)-1} R_{(n-k)\dvtx\varnothing}(1), &\quad $n > 0$,\vspace*{2pt}\cr
\displaystyle \prod_{k=0}^{|\bfi|-1} R_{n\dvtx\bfi|_{k}}(\bfi[k+1]), &\quad $n =
0$,
\vspace*{2pt}\cr
\displaystyle \prod_{k=0}^{|\bfi|-1} R_{n\dvtx\bfi|_{k}}(\bfi[k+1]) \prod_{k=0}^{|n|-1}
R_{-k\dvtx\varnothing}(1), &\quad $n < 0$.}
\]
Here we have used the convention that $\prod_{k=0}^{-1} x_k = 1$, to
deal with the case $|\bfi|=0$.

Let $\calW_{n\dvtx\bfi}$ be branching random walk limit associated with
crossing $n\dvtx\bfi$; see Lemma~\ref{BRWlem}. Then the duration of
crossing $n\dvtx\bfi$ is
%
\begin{equation}\label{Yxingeqn}
\calD_{n\dvtx\bfi} = \rho_{n\dvtx\bfi} \calW_{n\dvtx\bfi}.
\end{equation}
We summarize conditions for existence and continuity of $Y$ in the
theorem below.
%
\begin{theorem}\label{ThmMEBP1}
Suppose we are given subcrossing orientation distributions $p^\pm_{A}$
and weight distributions $F^\pm_{R|a}$, satisfying Assumptions \ref
{AssGW},~\ref{Ass1}\break and~\ref{assInf}.
Then there exists a continuous EBP process $Y$ with subcrossing
orientation distributions $p^\pm_{A}$ and crossing durations $\calD
_{n\dvtx\bfi} \stackrel{\mathit{fdd}}{=} \rho_{n\dvtx\bfi} \calW_{n\dvtx\bfi}$.

For each $n$, conditioned on the crossing orientations $\al^n_k$, the
random variables $\calW^n_k$ are mutually independent, and $\calW^n_k$
is conditionally independent of all $(A^m_j, R^m_j)$ for $m > n$.
Also, $\E(\calW^n_k  |  \al^n_k = i) = v^i$, and the distribution of
$\calW^n_k$ depends only on $\al^n_k$.

We call $Y$ the multifractal embedded branching process (MEBP) defined
by $p^\pm_A$ and $F^\pm_{R|a}$.
\end{theorem}

As a corollary of our construction we also obtain a novel
Galton--Watson cascade measure $\zeta$ on $\R_+$, constructed by
mapping the cascade measure $\nu$ from the boundary of the (doubly
infinite) tree to $\R_+$, using the measure $\bar{\nu}$ as a reference.
[Where $\bar{\nu}$ is defined on $\dUp$ by $\bar{\nu}(\dUp
_{n\dvtx\bfi}) =
\mu^{n} W_\bfi$.]

Mandelbrot, Fisher and Calvet~\cite{MFC97} described a class of
multifractal processes such that
\[
Y(at) \stackrel{\mathrm{fdd}}{=} M(a)Y(t)
\quad\mbox{and}\quad M(ab) \stackrel{d}{=}
M_1(a) M_2(b),
\]
where $M_1$ and $M_2$ are independent copies of $M$.
Write $A$ for $M^{-1}$, and then we can re-express the scaling rule for
$Y$ as
%
\begin{equation}\label{MFCmulteqn}
Y(A(a)t) \stackrel{\mathrm{fdd}}{=} a Y(t) \quad
\mbox{and}\quad A(ab) \stackrel{d}{=}
A_1(a)A_2(b),
\end{equation}
where $A_1$ and $A_2$ are independent copies of $A$.
When constructing our MEBP~$Y$, we noted that $Y^{(n)}(t)$ is
distributed as $2^nY^{(0)}(t \rho^{-n}_1)$.
More generally we have $Y^{(m+n)}(t) \stackrel{\mathrm{fdd}}{=} 2^nY^{(m)}(t
\rho
^{-n}_1)$, so sending $m\to\infty$ we get, for $n = 0, 1,
\ldots,$
\[
Y(t) \stackrel{\mathrm{fdd}}{=} 2^nY(t \rho^{-n}_1).
\]
This is close to the form (\ref{MFCmulteqn}) with $A(2^{-n}) = \rho
^{-n}_1 = \prod_{k=0}^{-n+1} R^k_1(1)$. The differences are that $A(a)$
is only defined for $a = 2^{-n}$, $n \in\Z_+$, and the product form
$A(ab) \stackrel{d}{=} A_1(a)A_2(b)$ does not quite hold because of the
dependence of $R^k_1$ on the orientation $\al^k_1$. [In fact, the
sequence $\{ (-\log\rho^{-n}_1, \al^{-n}_1) \}$ is Markov additive.]
Nonetheless, we recognize that MEBP processes possess a form of
discrete multifractal scaling. The full multifractal spectrum is
obtained in a forthcoming paper~\cite{DHJ}.

\section{On-line simulation}\label{secsim}

There are many ways we could make a multifractal time-change of a CEBP.
However, by defining the time-change via the crossing tree, we obtain a
fast on-line algorithm to simulate the process.
As before, we will suppose that we are given subcrossing orientation
distributions $p^\pm_A$ and weight distributions $F^\pm_{R|a}$,
satisfying Assumptions~\ref{AssGW},~\ref{Ass1} and~\ref{assInf}.
Let $Y$ be the corresponding MEBP. Then we will simulate the sequence
$\{(\calT^0_k, Y(\calT^0_k))\}$.
That is, we will simulate $Y$ at the \textit{spatial} scale of 1.
Given the multifractal nature of the process, the choice spatial scale
is not a restriction, as the process can be scaled to any desired resolution.
An immediate consequence of the definition of the crossing times $\calT
^0_k$ is the following bound on $Y$:
\[
Y(t) \in\bigl(Y(\calT^0_k)-1, Y(\calT^0_k)+1\bigr) \qquad\mbox{for } t \in(\calT
^0_k, \calT^0_{k+1}).
\]

The basis of our simulation is a Markov process, which describes the
line of descent of the current level zero crossing, from the spine down
to level 0.
For $n \geq m$ and $k \geq0$ let $\ka(m,n,k)$ be such that $C^m_k$ is
a subcrossing of $C^n_{\ka(m,n,k)}$, and let $S^n_k \in\{1, \ldots,
Z^{n+1}_{\ka(n,n+1,k)}\}$ be the position of $C^n_k$ within
$C^{n+1}_{\ka(n,n+1,k)}$.
Using this notation, if $n\dvtx\bfi$ is the tree-index of $C^0_k$, then
for $0 \leq m \leq n-1$, $\bfi[n-m] = S^m_{\ka(0,m,k)}$.
Let $\calY^n(k) = (\ka(0,n,k), S^n_{\ka(0,n,k)}, Z^{n+1}_{\ka
(0,n+1,k)}, A^{n+1}_{\ka(0,n+1,k)},\break R^{n+1}_{\ka(0,n+1,k)})$, which is
a description of the level $n$ super-crossing of $C^0_k$, and the
family it belongs to.

Let $N(k)$ be the smallest $n$ such that $\ka(0,n+1,k) = 1$, and put
\[
\calY(k) = \bigl(\calY^0(k), \ldots, \calY^{N(k)}(k)\bigr).
\]

\begin{lem}
$\calY$ is a Markov process.
\end{lem}
\begin{pf}
We first show how to update $\calY(k)$ to obtain $\calY(k+1)$.
Let $M$ be the largest $m \leq N(k)$ such that
\[
S^n_{\ka(0,n,k)} = Z^{n+1}_{\ka(0,n+1,k)} \qquad\mbox{for } n = 0, \ldots
, m.
\]
That is, for all $m \leq M$ we have that $C^m_{\ka(0,m,k)}$ is the last
level $m$ crossing in its family.

If $M = N(k)$, then $N(k+1) = N(k) + 1$, and $\calY$ gains the
component $\calY^{N(k+1)}(k+1)$.
Let $n = N(k+1)$. Then we have $\ka(0,n,k+1) = 2$, $S^{n}_2 = 2$ and
$\ka(0,n+1,k+1) = 1$.
The distribution of $(Z^{n+1}_1, A^{n+1}_1, R^{n+1}_1)$ depends on
$\calY(k)$ only through $\al^{n}_1 = A^{n+1}_1(1)$, which is given by
$A^{n}_1(Z^{n}_1)$.

If $M < N(k)$, then for $n = M + 1$ we have $\ka(0,n,k+1) = \ka(0,n,k)
+ 1$, $S^n_{\ka(0,n,k+1)} = S^n_{\ka(0,n,k)} + 1$ and $\ka
(0,n+1,k+1) =
\ka(0,n+1,k)$.
Thus $Z^{n+1}_{\ka(0,n+1,k+1)} = Z^{n+1}_{\ka(0,n+1,k)}$,
$A^{n+1}_{\ka
(0,n+1,k+1)} = A^{n+1}_{\ka(0,n+1,k)}$, and $R^{n+1}_{\ka(0,n+1,k+1)} =
R^{n+1}_{\ka(0,n+1,k)}$.\vspace*{1pt}

For $n > M+1$, we have $\calY^n(k+1) = \calY^n(k)$.

For $n = M, \ldots, 0$, we generate $\calY^n(k+1)$ recursively.
We have $\ka(0,n,k+1) = \ka(0,n,k)+1$, $S^{n}_{\ka(0,n,k+1)}=1$, and
$\ka(0,n+1,k+1) = \ka(0,n+1,k)+1$.
The distribution of $(Z^{n+1}_{\ka(0,n+1,k+1)}, A^{n+1}_{\ka
(0,n+1,k+1)}, R^{n+1}_{\ka(0,n+1,k+1)})$ is determined by $\al
^{n+1}_{\ka(0,n+1,k+1)}$, that is, $A^{n+2}_{\ka
(0,n+2,k+1)}(S^{n+1}_{\ka(0,n+1,k+1)})$.
Thus $\{ \calY^n(k+1) \}_{n=0}^M$ depends on $\calY(k)$ only through
$\al^{M+1}_{\ka(0,M+1,k+1)} = A^{M+2}_{\ka(0,M+2,k+1)}(S^{M+1}_{\ka
(0,M+1,k+1)}) =
A^{M+2}_{\ka(0,M+2,k)}(S^{M+1}_{\ka(0,M+1,k)+1})$.\vspace*{2pt}

That $\calY$ is Markov follows from the conditional independence of the
$(Z^n_k, A^n_k$, $R^n_k)$ given the orientations $\al^n_k$.
\end{pf}

From $\calY(k)$, we get the orientation $\al^0_k$ of $C^0_k$, and the
weights
\[
R^{n+1}_{\ka(0,n+1,k)}\bigl(S^n_{\ka(0,n,k)}\bigr)\qquad \mbox{for $n=0, \ldots
, N(k)$}.
\]
To calculate the crossing duration $\calD^0_k$ we also need
$R^{n+1}_1(1)$, for $n=0, \ldots, \break N(k)$ and $\calW^0_k$.
Keeping track of the spine weights $R^{n+1}_1(1)$ is no problem.
Calculating $\calW^0_k$ is less straightforward.
We do have that the $\calW^0_k$ are conditionally independent given the
$\al^0_k$, but we do not have an explicit formulation of the density of
$(\calW^0_k  |  \al^0_k=\pm)$.

The simplest way to approximate the $\calW^0_k$ is to generate a BRW
(using $p^\pm_A$ and $F^\pm_{R|a}$) for a fixed number of generations,
$m$ say, and sum the node weights across the final generation.
However, this is exactly the same as setting the $\calW^0_k$ to be
constant, then scaling the resulting process by $2^{-m}$, so we will
just set $\calW^0_k$ equal to its mean $v^{\al^0_k}$.
%
\begin{remark}
Writing $Y$ as $X \circ\calM^{-1}$, where $X$ is the CEBP
corresponding to $Y$, we note that $X$ and $\calM$ are, in general, dependent.
However, in the case where $X$ is Brownian motion, we can construct
$\calM$ independently of $X$, simply by taking the orientations $\al
^0_k$ as i.i.d. random variables, equal to $+$ and $-$ with equal probability.
This is because for Brownian motion $X(T^0_k)$ is just a simple random walk.
In fact, in this case, there need not be any relation at all between
the crossing tree of $X$ and that used to construct~$\calM$.
\end{remark}

\subsection{Pseudo-code}

We give pseudo code for simulating $\{(\calT^0_k, Y(\calT^0_k))\}$,
with the crossing durations $\calD_{n\dvtx\bfi}$ approximated by $\rho
_{n\dvtx\bfi}\E(\calW_{n\dvtx\bfi})$ (i.e., putting $\calW_{n\dvtx\bfi} = v^i$,
where $i = \al_{n\dvtx\bfi}$).

\begin{figure}

\includegraphics{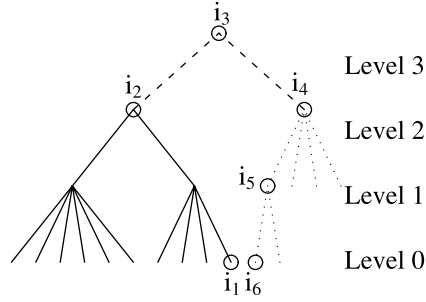}

\caption
{Action of the procedures \texttt{Increment} and \texttt{Expand}.
Suppose that at iteration $k$ we have generated the tree given by solid
black lines only, so that $N(k)=1$, and $\calY^{N(k)}(k)$ describes the
family of node $\bfi_1$.
When we reach node $\bfi_1$, we are at the end of level 0 and 1 crossings.
To generate the next level 0 node, we first need to increase $N(k)$ by
1 and generate the family of node $\bfi_3$, which is the role of the
procedure \texttt{Expand}.
Next, procedure \texttt{Increment} goes down the tree and generates the
families of nodes $\bfi_4$ and $\bfi_5$, hence generating $\bfi_6$.}
\label{figalgo}
\end{figure}

Updating $\calY(k)$ is handled by procedures \texttt{Expand} and
\texttt
{Increment}.
Procedure \texttt{Expand} checks if $M = N(k)$.
If so, it then generates the component $\calY^{M+1}(k)$ and updates $N(k)$.
Assuming $M < N(k)$, procedure \texttt{Increment} updates $\calY^n(k)$
to $\calY^n(k+1)$ recursively, for $n = M+1, \ldots, 0$.
The actions of \texttt{Expand} and \texttt{Simulate} are illustrated in
Figure~\ref{figalgo}.

Given sample position $Y(\calT^0_{k})$, sample time $\calT^0_{k}$ and
crossing state $\calY(k)$, the procedure \texttt{Simulate} applies the
procedures \texttt{Expand} and \texttt{Increment}, calculates
$Y(\calT
^0_{k+1})$, $\calT^0_{k+1}$ and $\calY(k+1)$, then increments $k$.
Procedure \texttt{Initialize} generates an initial $Y(\calT^0_{1})$,
$\calT^0_{1}$ and $\calY(1)$ suitable for passing to \texttt
{Simulate}.

Recall that $u = \Pb( \al^{n+1}_1 = +  |  \al^n_1 = +)$ and $v =
\Pb(
\al^{n+1}_1 = +  |  \al^n_1 = -)$.
Here $\al^{N(k)+1}_1$ is given by
$A^{N(k)+1}_1(Z^{N(k)+1}_1)$.\vspace*{12pt}

\texttt{Procedure Expand} $\calY(k)$\vspace*{2pt}

\hspace{1cm} \texttt{If} $S^{N(k)}_{\ka(0,N(k),k)} = Z^{N(k)+1}_{\ka(0,N(k)+1,k)}$
\texttt{Then}\vspace*{2pt}

\hspace{1.5cm} $\kappa(0,N(k)+2,k) = 1$\vspace*{2pt}

\hspace{1.5cm} Generate $\al^{N(k)+2}_1$ using $u$, $v$ and
$\al^{N(k)+1}_1$\vspace*{2pt}

\hspace{1.5cm} Generate $(Z^{N(k)+2}_{\ka(0,N(k)+2,k)}, A^{N(k)+2}_{\ka(0,N(k)+2,k)},
R^{N(k)+2}_{\ka(0,N(k)+2,k)})$\vspace*{2pt}

\hspace{2cm} using the distributions $p^i_A$ and $F^i_{R|a}$\vspace*{1pt}

\hspace{2cm} conditioned on the first offspring having orientation
$\al^{N(k)+1}_1$\vspace*{1pt}

\hspace{2cm} where $i = \al^{N(k)+2}_1 \in \{+, -\}$\vspace*{2pt}

\hspace{1.5cm} $S^{N(k)+1}_{\ka(0,N(k)+1,k)} = 1$\vspace*{2pt}

\hspace{1.5cm} Store $R^{N(k)+2}_1(1)$\vspace*{1pt}

\hspace{1.5cm} $N(k) = N(k) + 1$\vspace*{1pt}

\hspace{1cm} \texttt{End If}

\texttt{End Procedure}\vspace*{12pt}

\texttt{Procedure Increment} $\calY^n(k)$\vspace*{1pt}

\hspace{1cm} \# Assume that $C^{n-1}_k$ is at the end of a level $n$
crossing,\vspace*{1pt}

\hspace{1cm} \# so $S^{n-1}_{\ka(0,n-1,k)} = Z^n_{\ka(0,n,k)}$. This is always the case for
$n=0$\vspace*{2pt}

\hspace{1cm} $\kappa(0,n,k+1)=\kappa(0,n,k)+1$\vspace*{2pt}\vadjust{\goodbreak}

\hspace{1cm} \texttt{If} $S_{\kappa(0,n,k)}^{n}=Z_{\kappa(0,n+1,k)}^{n+1}$
\texttt{Then}\vspace*{2pt}

\hspace{1.5cm} \texttt{Increment} ${\cal X}^{n+1}(k)$\vspace*{1pt}

\hspace{1.5cm} $S^{n}_{\kappa(0,n,k+1)}=1$\vspace*{1pt}

\hspace{1.5cm} Generate $(Z^{n+1}_{\kappa(0,n+1,k+1)}, A_{\kappa(0,n+1,k+1)}^{n+1},
R^{n+1}_{\kappa(0,n+1,k+1)})$\vspace*{2pt}

\hspace{2cm} using the distributions $p^i_A$ and $F^i_{R|a}$\vspace*{1pt}

\hspace{2cm} where $i = A_{\kappa(0,n+2,k+1)}^{n+2}(S^{n+1}_{\kappa(0,n+1,k+1)}) \in
\{+,-\}$\vspace*{1pt}

\hspace{1cm} \texttt{Else}\vspace*{1pt}

\hspace{1.5cm} ${\cal X}^q(k+1)=\calY^q(k)$ \texttt{for}
$q=n+1,\ldots,N(k)$\vspace*{2pt}

\hspace{1.5cm} $S^{n}_{\kappa
(0,n,k+1)}=S^{n}_{\kappa(0,n,k)}+1$\vspace*{2pt}

\hspace{1cm} \texttt{End If}

\texttt{End Procedure}\vspace*{12pt}

We apply procedure \texttt{Increment} to $\calY^0(k)$, and then it is
recursively applied to all $\calY^n(k)$ such that $C^q_{\ka(0,q,k)}$ is
at the end of a level $q+1$ crossing for all $0 \leq q < n$.
$\calY^n(k+1) = \calY^n(k)$ for all $n$ larger than this.\vspace*{12pt}

\texttt{Procedure Simulate}

\hspace{1cm} \texttt{Expand} $\calY(k)$\vspace*{1pt}

\hspace{1cm} \texttt{Increment} ${\cal X}^0(k)$\vspace*{1pt}

\hspace{1cm} Put $i = A^1_{\ka(0,1,k+1)}(S^0_{k+1})$\vspace*{2pt}

\hspace{1cm} \texttt{If} $i = +$ \texttt{Then}\vspace*{1pt}

\hspace{1.5cm} $Y(\calT^{0}_{k+1}) = Y(\calT^{0}_{k}) + 1$\vspace*{2pt}

\hspace{1cm} \texttt{Else}\vspace*{1pt}

\hspace{1.5cm} $Y(\calT^{0}_{k+1}) = Y(\calT^{0}_{k}) - 1$\vspace*{2pt}

\hspace{1cm} \texttt{End If}\vspace*{1pt}

\hspace{1cm} $\calT^0_{k+1} = \calT^0_k + v^i \prod_{j=0}^{N(k+1)} ( R^{j+1}_{\kappa(0,j+1,k+1)}( S^j_{\kappa(0,j,k+1)}) /
R^{j+1}_1(1))$\vspace*{2pt}

\hspace{1cm} $k \leftarrow k + 1$\vspace*{1pt}

\texttt{End Procedure}\vspace*{12pt}

To initialize the algorithm, the procedure \texttt{Initialize} is used.
Recall that $(v^+, v^-)^T$ is the right $\mu(1)$-eigenvector of
$M(1)$.\vadjust{\goodbreak}\vspace*{12pt}

\texttt{Procedure Initialize} $\calY(1)$\vspace*{1pt}

\hspace{1cm} $k = 1$, $N(1) = 0$, $\kappa(0,0,1)=1$,
$\kappa(0,1,1)=1$\vspace*{1pt}

\hspace{1cm} Put $\al^1_1 = i = +$ with probability $a$\vspace*{1pt}

\hspace{1cm} Generate $(Z^1_1, A^1_1, R^1_1)$ using the
distributions\vspace*{1pt}

\hspace{1.5cm} $p^i_A$ and $F^i_{R|a}$, with $i = \al^1_1$\vadjust{\goodbreak}

\hspace{1cm} $S^0_1 = 1$\vspace*{1pt}

\hspace{1cm} Store $R^1_1(1)$\vspace*{1pt}

\hspace{1cm} $\calT^0_1 = v^i$\vspace*{1pt}

\hspace{1cm} \texttt{If} $i = +$ \texttt{Then} $Y(\calT^0_1) = 1$ \texttt{Else} $Y(\calT^0_1) = -1$ \texttt{End
If}\vspace*{1pt}

\texttt{End Procedure}\vspace*{12pt}

\begin{figure}[b]

\includegraphics{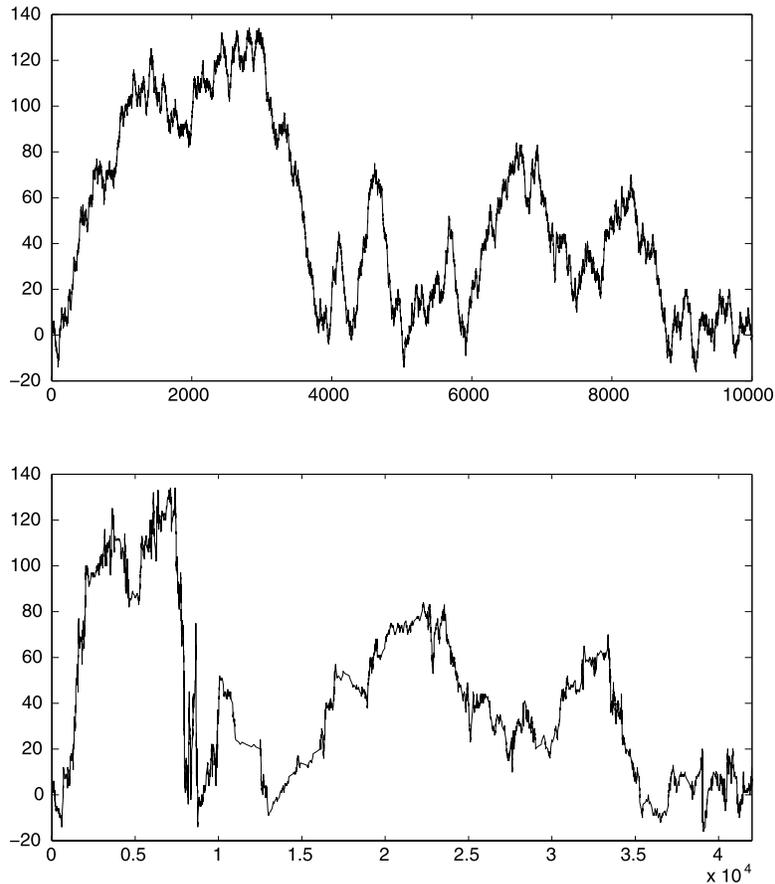}

\caption
{Top figure: CEBP process where the offspring consist of a
\textup{geometric}($0.6$) number of excursions, each up--down or down--up with
equal probability, followed by either an up--up or down--down direct
crossing [compare this with Brownian motion, for which there are a
\textup{geometric}($0.5$) number of excursions].
Bottom figure: MEBP process obtained from a multifractal time change of
the top CEBP process, with i.i.d. gamma distributed weights.}
\label{EBPMEBP}
\end{figure}

An implementation is available from the web page of Jones~\cite{J}.
An example of the type of signal obtained with this algorithm is given
in Figure~\ref{EBPMEBP}, where we have represented an MEBP process
%
with its corresponding CEBP.
$p^\pm_A$ and $F^\pm_{R|a}$ are described in the caption.

\subsection{Efficiency}

Consider the tree descending from crossing $C^{N(k)}_1$ down to level 0.
On average $C^{N(k)}_1$ has $\mu^{N(k)}$ level 0 subcrossings, so we
must have $N(k) = O(\log k)$.
At each step, the number of operations required by procedure \texttt
{Expand} is fixed [independent of $N(k)$], but we can go through
\texttt
{Increment} up to $N(k)$ times, so the number of operations required by
\texttt{Simulation} is of order $N(k)$.
Thus, to generate $n$ steps, we use $O(n \log n)$ operations, since
$\sum _{k=1}^n \log k = O(n\log n)$, and $O(\log n)$ storage.
The algorithm is on-line, meaning that given the current state [of size
$O(\log n)$] we can generate the next immediately [using $O(\log n)$
operations].

\section{Randomizing the starting point}

Crossing times are points where the behavior of the process can change,
spatially and temporally, and the higher the level, the more dramatic
this can be.
For MEBP processes, 0 is a crossing time for all levels, and because of
this we cannot expect MEBP to have stationary increments.
To avoid the problem of 0 being special, we would like to start the
process at a ``random'' time, as if the process had been running since
time immemorial and we just happened across it.

To make the idea of a ``random'' starting time more precise, let $Y$ be
an MEBP and $\{ Y^{(n)} \}$ the nested sequence of processes used to
construct~$Y$, where $Y^{(n)}$ is a single level $n$ crossing from $0$
to $\pm2^n$.
Choose a time $t$ uniformly in $[0, T^n_1] = [0, \calD_{n\dvtx\varnothing}]$.
For any $\bfi\in\Up_{n\dvtx\varnothing}$, the probability that $t$ is in
$C_{n\dvtx\bfi}$ is proportional to the crossing duration $\calD_{n\dvtx\bfi
} =
\rho_{n\dvtx\bfi} \calW_{n\dvtx\bfi}$.
That is, choosing $t$ is equivalent to choosing $n\dvtx\bfj\in\dUp
_{n\dvtx\varnothing}$ so that the probability that $n\dvtx\bfj|_{|\bfi|} =
n\dvtx\bfi$ is proportional to $\rho_{n\dvtx\bfi} \calW_{n\dvtx\bfi}$.
It turns out that we can do exactly this using a size-biased measure
for a multitype branching random walk.

Size-biased measures for branching processes were introduced by Lyons,
Pemantle and Peres~\cite{LPP95} and generalized to branching random
walks by Lyons~\cite{L97}.
Kyprianou and Sani~\cite{KS01} then extended their construction to
multitype branching random walks.
Fix $n$, and for brevity write $\bfi$ for $n\dvtx\bfi$.
Let $\Om$ be the space of marked trees, where the mark associated with
node $\bfi$ is $(-\log R_{\bfi|_{k-1}}(\bfi[k]), \al_\bfi)$, writing
$k$ for~$|\bfi|$.
Let $\calF$ be the $\sigma$-field generated by all finite truncations
of trees.
The offspring orientation distributions $p^\pm_A$ and weight
distributions $F^\pm_{R|a}$ induce a measure $\xi$ on $(\Om, \calF)$.
Let $\tilde{\Om}$ be the space of trees with a distinguished line of
descent $\bfi\in\dUp$, called a spine, and $\tilde{\calF}$ the
$\sigma
$-field generated by all finite truncations of trees with spines.
Kyprianou and Sani define a size-biased measure $\tilde{\pi}$ on
$(\tilde{\Om}, \tilde{\calF})$ such that
%
\begin{equation}\label{eqnsize-biased}
\int_{\bfj\in\dUp_\bfi} d\tilde{\pi}(\Up, \bfj) = \frac{\rho
_\bfi\calW
_\bfi}{v^{\al_\varnothing}} \,d\xi(\Up).
\end{equation}
This is precisely what we want, and, remarkably, the measure can be
constructed using the original multitype branching walk, modified so
that the offspring generation down the spine is size-biased.
That is, rather than construct $Y^{(n)}$ and then choose a spine, we
can construct the process and the spine together.

Let $\bfx\in\dUp$ be the spine, and let $\tilde{p}^\pm_{A}$ and
$\tilde{F}^\pm_{R|a}$ be the offspring orientation and weight
distributions for nodes on the spine. Then from~\cite{KS01}, Section 2,
we have that
\begin{eqnarray*}
\tilde{p}^i_{A}(a) \tilde{F}^i_{R|a}(r)
&=& \Pb(A_{\bfx|_n} = a, R_{\bfx|_n} \leq r  |  \al_{\bfx|_n} =
i) \\
&\propto& p^i_{A}(a) \sum_{j=1}^{|a|} v^{a(j)} \int_{s \leq r} s(j)
F^i_{R|a}(ds).
\end{eqnarray*}
Note here that $s$ and $r$ are in $\R_+^{|a|}$.
Putting $r = \infty^{|a|}$ to get $\tilde{p}^i_{A}(a)$, and then
dividing out $\tilde{p}^i_{A}(a)$ to get $\tilde{F}^i_{R|a}(r)$,
gives us
\begin{eqnarray*}
\tilde{p}^i_{A}(a) &\propto& p^i_A(a) \sum_{j=1}^{|a|} v^{a(j)} \int
_{\R
_+^{|a|}} s(j) F^i_{R|a}(ds), \\
\tilde{F}^i_{R|a}(r) &\propto& \sum_{j=1}^{|a|} v^{a(j)} \int_{s
\leq
r} s(j) F^i_{R|a}(ds).
\end{eqnarray*}
That these are well defined follows from Assumption~\ref{Ass1}.

In the case where the offspring weights are i.i.d. with distribution
$F$, we get
\begin{eqnarray*}
\tilde{p}^i_{A}(a) &\propto& |a| p^i_A(a), \\
\tilde{F}^i_{R|a}(r) &\propto& \sum_{j=1}^{|a|} \int_0^{r(j)} s F(ds)
\prod_{i\neq j} F(r(i)).
\end{eqnarray*}
The first of these is clearly a size-biased version of $p^i_A$.
The second can be interpreted as conditioning on which offspring is on
the spine, then size-biasing the weight for that offspring.

For selecting the next node on the spine, we again have from \cite
{KS01}, Section~2, that
\[
\tilde{p}_{a,r}(j):= \Pb(\bfx[n+1] = j  |  A_{\bfx|_n}=a,
R_{\bfx
|_n}=r) \propto v^{a(j)} r(j).
\]

Kyprianou and Sani also also show that under $\tilde{\pi}$, the
sequence $\{ \al_{\bfx|_n} \}_{n=1}^\infty$ of orientations down the
spine is Markovian, with transition probabilities
\[
\pmatrix{v^+&0\cr
0&v^-}^{-1}
M(1)
\pmatrix{v^+&0\cr
0&v^-}.
\]
The stationary distribution is $(u^+v^+, u^-v^-)$, and so the reversed
chain (moving up the spine) has transition matrix
%
\begin{equation}\label{eqnspineorientation}
\pmatrix{u^+&0\cr
0&u^-}
^{-1}
M(1)^T
\pmatrix{u^+&0\cr
0&u^-},
\end{equation}
and the same stationary distribution as before.
Note that it follows from assumptions~\ref{AssGW} and~\ref{Ass1} that
$\bfu, \bfv> 0$.

\subsection{MEBP construction with random start}

We now show how, given an MEBP $Y\dvtx [0,\infty) \to\R$ generated by
$p^\pm_A$ and $F^\pm_{R|a}$, we can construct a shifted version,
$\tilde
{Y}\dvtx (-\infty, \infty) \to\R$, with a ``randomly'' chosen starting point.
Where unambiguous, we will use the same notation to describe $\tilde
{Y}$ as $Y$, and we will assume that assumptions~\ref{AssGW} and \ref
{Ass1} hold throughout.
As before, we start by constructing a crossing of size 1 (level 0).
Let $\bfx$ be the spine, which will be the line of descent
corresponding to time 0.
Accordingly, we will write $C^{-n}_0 = C_{\bfx|_n}$ for the level $-n$
spinal crossing.
Note that previously, the first crossing at level $-n$ was labeled 1,
and started at time 0.
For our new construction, time 0 will occur somewhere in the interior
of crossing $C^{-n}_0$, so crossing $C^{-n}_1$ will still be the first
full crossing to occur after time 0.

The generation $n$ (level $-n$) nodes in $\Up_n$ are totally ordered
according to the rule $\bfi< \bfj$ if and only if, for some $m$,
$\bfi
|_m = \bfj|_m$ and $\bfi[m+1] < \bfj[m+1]$.
For $\bfi, \bfj\in\Up_n$ let
\[
d(\bfi, \bfj) = \cases{
|\{ \bfk\dvtx  \bfi< \bfk\leq\bfj\}|, &\quad $\bfi< \bfj$,\cr
0, &\quad $\bfi= \bfj$,\cr
-|\{ \bfk\dvtx  \bfi> \bfk\geq\bfj\}|, &\quad $\bfi> \bfj$.}
\]
We will write $C^{-n}_{d(\bfx,\bfi)}$ for $C_{\bfi}$.

Set the orientation of $C^0_0$ to be $+$ with probability $u^+v^+$, and
then generate $(A^0_0, R^0_0)$ using $\tilde{p}^i_A$ and $\tilde
{F}^i_{R|a}$, where $i = \al^0_0$.
Choose $j \in\{1, \ldots, Z_\varnothing\}$ using $\tilde
{p}_{A_\varnothing
, R_\varnothing}$, and then put $\bfx|_1 = j$.
Subsequent generations are produced using $p^\pm_A$ and $F^\pm_{R|a}$
for nodes off the spine, and $\tilde{p}^\pm_A$ and $\tilde{F}^\pm
_{R|a}$ for the spinal node.
The spinal node in the next generation is chosen using $\tilde{p}_{a,r}$.
Crossing durations are defined as before; that is, $\calD^{-n}_k =
\rho
^{-n}_k \calW^{-n}_k$, where $\calW_\bfi$ is the $\tilde{\pi
}$-a.s.
limit of $\sum_{\bfj\in\Up_n\cap\Up_\bfi} \rho_\bfj/\rho_\bfi$.
For $k \neq0$ (nodes off the spine) the convergence of this sequence
a.s. and in mean follows as before.
For $k = 0$ (nodes on the spine) a.s. convergence follows from
(\ref{eqnsize-biased}) and the fact that $\rho_{\bfx|_n}
\calW
_{\bfx|_n} \in(0, \infty)$ $\xi$-a.s.

Given crossing durations, we define crossing times as follows.
Time 0 corresponds to the spine $\bfx$.
For any $m \geq0$, $\calT^{-m}_0 > 0$ is the first time the process
starts a level $-m$ crossing:
\begin{eqnarray*}
\calT^{-m}_0
&=& \lim_{n\to\infty} \sum_{\bfi\in\Up_n, \bfi|_m=\bfx|_m,
\bfi> \bfx}
\rho_\bfi\calW_\bfi,\\
\calT^{-m}_{k+1} &=& \calT^{-m}_k + \rho^{-m}_{k+1} \calW^{-m}_{k+1}
\qquad\mbox{for } k\geq0 ,\\[-2pt]
\calT^{-m}_{-1}
&=& \lim_{n\to\infty} \sum_{\bfi\in\Up_n, \bfi|_m=\bfx|_m,
\bfi< \bfx}
\rho_\bfi\calW_\bfi,\\[-2pt]
\calT^{-m}_{-k-1} &=& \calT^{-m}_{-k} - \rho^{-m}_{-k} \calW^{-m}_{-k}
\qquad\mbox{for } k\geq1.
\end{eqnarray*}
We also put $\tilde{Y}(0) = 0$ and
\begin{eqnarray*}
\tilde{Y}(\calT^{-m}_0)
&=& \lim_{n\to\infty} \sum_{\bfi\in\Up_n, \bfi|_m=\bfx|_m,
\bfi> \bfx}
\al_\bfi2^{-n}, \\[-2pt]
\tilde{Y}(\calT^{-m}_{k+1}) &=& \tilde{Y}(\calT^{-m}_k) + \al
^{-m}_{k+1} 2^{-m} \qquad\mbox{for } k\geq0, \\[-2pt]
\tilde{Y}(\calT^{-m}_{-1})
&=& \lim_{n\to\infty} \sum_{\bfi\in\Up_n, \bfi|_m=\bfx|_m,
\bfi< \bfx}
\al_\bfi2^{-n}, \\[-2pt]
\tilde{Y}(\calT^{-m}_{-k-1}) &=& \tilde{Y}(\calT^{-m}_{-k}) + \al
^{-m}_{-k} 2^{-m} \qquad\mbox{for } k\geq1.
\end{eqnarray*}
So for $k\geq1$, $C^{-m}_k$ is from $\calT^{-m}_k$ to $\calT
^{-m}_{k+1}$, while for $k\leq0$ it is from $\calT^{-m}_{k-1}$ to~$\calT^{-m}_k$.

Let $\tilde{Y}^{(0)}$ be the level 0 crossing constructed above. We now
show how to extend the construction from $\tilde{Y}^{(n)}$ to
$\tilde{Y}^{(n+1)}$. Let $n\dvtx\bfx$ be the spine starting at level
$n$. First choose $\al^{n+1}_0 = i$ using the reversed Markov chain
\ref {eqnspineorientation}, then choose $(A^{n+1}_0, R^{n+1}_0)$ and
$(n+1)\dvtx\bfx[1] = j$ using $\tilde{p}^i_A$, $\tilde{F}^i_{R|a}$ and
$\tilde{p}_{a,r}$, all conditioned on $\al^n_0$, which is the
orientation of $(n+1)\dvtx\bfx[1]$. Put the $j$th level $n$ subcrossing
of $\tilde{Y}^{(n+1)}$, that is~$C^n_0$, equal to $\tilde{Y}^{(n)}$.
For the other level $n$ subcrossings, we use the construction of
Section~\ref{secextend}, and scale the $k$th subcrossing by
$R^{n+1}_0(k)/R^{n+1}_0(j)$. That is, we use the weights up the spine,
from level 0 to $n$, to rescale the process. Let $\tilde{Y}$ be the
limit of the $\tilde{Y}^{(n)}$.

To see that $\tilde{Y}(t)$ is defined for all $t \in\R$ we need two things.
First we note that from the form of $\tilde{p}_{a,r}$, with probability
1 we cannot have $n\dvtx\bfx[1]$ equal to 1 eventually, or equal to
$Z_{n\dvtx\varnothing}$ eventually.
That is, at all levels there will be crossings to the left and right of
the spinal crossing.
Second, we need to know that the scaling coming from the spine weights
grows to infinity, that is, $\prod_{k=1}^{n+1} R^k_0(k\dvtx\bfx[1]) \to
0$ a.s. as $n \to\infty$.

As noted above, the sequence of orientations up the spine is a Markov process.
Because the weights are conditionally independent given the
orientations, the sequence $(\sum_{k=1}^{n+1} \log R^k_0(k\dvtx\bfx[1]),
\al^{n+1}_0)$ is Markov additive.
Thus, $\sum_{k=1}^{n+1} \log R^k_0(k\dvtx\bfx[1]) \to-\infty$ a.s.,
equivalently $\prod_{k=1}^{n+1} R^k_0(k\dvtx\bfx[1]) \to0$ a.s.,
provided the expected increments of the sum are negative.
That is, provided the following assumption holds (this replaces
Assumption~\ref{assInf}).
%
\begin{assp}\label{assspinalweights}
Let $R^\pm$ be a random spinal weight, chosen according to $\tilde
{p}^\pm_A$, $\tilde{F}^\pm_{R|a}$ and $\tilde{p}_{a,r}$.
Then we assume that
\[
u^+v^+ \E\log R^+ + u^-v^- \E\log R^- < 0.\vadjust{\goodbreak}
\]
\end{assp}

It remains an open problem to show that the process $\tilde{Y}$ has
stationary increments.\vspace*{-2pt}

\subsection{On-line simulation}

To simulate $\tilde{Y}$ we need only modify procedures \texttt{Expand}
and \texttt{Initialie}.
Note that the spinal crossings are now counted as crossing 0 at each
level, so $N(k)$ is the smallest $n$ such that $\ka(0,n+1,k) =
0$.\vspace*{11pt}

\texttt{Procedure Expand} $\calY(k)$\vspace*{2pt}

\hspace{1cm} \texttt{While} $S^{N(k)}_{\ka(0,N(k),k)} =
Z^{N(k)+1}_{\ka
(0,N(k)+1,k)}$ \texttt{Do}\vspace*{2pt}

\hspace{1.5cm} $\kappa(0,N(k)+2,k) = 0$\vspace*{2pt}

\hspace{1.5cm} Generate $\al^{N(k)+2}_0$ using $(u^+v^+, u^-v^-)$ and
$\al^{N(k)+1}_0$\vspace*{2pt}

\hspace{1.5cm} Generate $A^{N(k)+2}_0$, $R^{N(k)+2}_0$ and
$S^{N(k)+1}_0$\vspace*{2pt}

\hspace{2cm} using the distributions $\tilde{p}^i_A$, $\tilde
{F}^i_{R|a}$ and $\tilde{p}_{a,r}$\vspace*{2pt}

\hspace{2cm} conditioned on offspring $S^{N(k)+1}_0$ having orientation
$\al^{N(k)+1}_0$\vspace*{2pt}

\hspace{2cm} where $i = \al^{N(k)+2}_0 \in\{+, -\}$\vspace*{2pt}

\hspace{1.5cm} Store $R^{N(k)+2}_0(S^{N(k)+1}_0)$\vspace*{2pt}

\hspace{1.5cm} $N(k) = N(k) + 1$\vspace*{2pt}

\hspace{1cm} \texttt{End While}

\texttt{End Procedure}

\texttt{Procedure Initialize} $\calY(0)$

\hspace{1cm} $k = 0$, $N(0) = 0$, $\kappa(0,0,0)=0$,
$\kappa(0,1,0)=0$\vspace*{2pt}

\hspace{1cm} Put $\al^1_0 = +$ with probability $u^+v^+$\vspace*{2pt}

\hspace{1cm} Generate $A^1_0$, $R^1_0$ and $S^0_0$ using the
distributions\vspace*{2pt}

\hspace{1.5cm} $\tilde{p}^i_A$, $\tilde{F}^i_{R|a}$ and $\tilde
{p}_{a,r}$, with $i = \al^1_0$\vspace*{2pt}

\hspace{1cm} Store $R^1_0(S^0_0)$\vspace*{2pt}

\hspace{1cm} $\calT^0_0 = 0$, $Y(\calT^0_0) = 0$\vspace*{2pt}

\texttt{End Procedure}\vspace*{-2pt}

\section*{Acknowledgments}

The authors are grateful for the many constructive comments received
from their anonymous referees.\vspace*{-2pt}



\printaddresses

\end{document}